\documentclass[a4paper,reqno]{amsart}
\usepackage[T1]{fontenc}
\usepackage[latin1]{inputenc}
\usepackage{graphicx}
\usepackage{amssymb}
\usepackage{latexsym}
\usepackage{mathrsfs}
\usepackage[all]{xy}
\newtheorem{thm}{Theorem}[section]
\newtheorem{cor}[thm]{Corollary}
\newtheorem{lem}[thm]{Lemma}
\newtheorem{prop}[thm]{Proposition}
\theoremstyle{definition}
\newtheorem{ex}[thm]{Example}
\theoremstyle{definition}
\newtheorem{defn}[thm]{Definition}
\theoremstyle{remark}
\newtheorem{rem}[thm]{Remark}
\numberwithin{equation}{section}
\newcommand{\norm}[1]{\left\Vert#1\right\Vert}
\newcommand{\abs}[1]{\left\vert#1\right\vert}

\newcommand{\Real}{\mathbb R}

\newcommand{\C}{\mathbb C}
\newcommand{\Z}{\mathbb Z}
\newcommand{\N}{\mathbb N}

\newcommand{\Sc}{\mathfrak{S}}
\newcommand{\Proj}{\mathbb P}

\newcommand{\To}{\longrightarrow}

\newcommand{\A}{\mathcal{A}}
\newcommand{\Al}{\mathscr{A}}

\newcommand{\Ev}{\mathscr{E}}
\newcommand{\Dom}{\mathscr{D}}
\newcommand{\U}{\mathcal{U}}
\newcommand{\F}{\mathcal{F}}
\newcommand{\Ge}{\mathcal{G}_e}
\newcommand{\Gr}{\mathrm{Gr}}

\newcommand{\GL}{\mathrm{GL}}
\newcommand{\Hilb}{\mathcal{H}}
\newcommand{\Ltwo}{\mathrm{L}^2}
\newcommand{\Lp}{\mathcal{L}}
\newcommand{\tensor}{\otimes}
\newcommand{\Lie}{\mathfrak{g}}
\newcommand{\lie}{\textrm{Lie}}
\newcommand{\Dirac}{\rlap{$D$}{\,/}}
\newcommand{\Cinfty}{\textrm{C}^\infty}

\newcommand{\Aut}{\textrm{Aut}}
\newcommand{\tr}{\mathrm{Tr}}
\renewcommand{\abs}[1]{\vert  #1 \vert}
\newcommand{\arr}{\rightrightarrows}
\newcommand{\id}{\mathrm{id}}
\newcommand{\Det}{\mathrm{Det}}
\newcommand{\fa}{\partial}
\newcommand{\de}{\sigma}
\newcommand{\Ob}{\mathrm{Ob}}
\newcommand{\Mat}[4]{
\left(
\begin{array}{ccc}
#1 & #2 \\ #3 & #4
\end{array}
\right)
}

\newcommand{\ket}[1]{\vert #1\rangle}
\newcommand{\stack}{\mathfrak{X}}
\newcommand{\gerbe}{\mathfrak{R}}
\newcommand{\Hom}{\textrm{Hom}}
\newcommand{\End}{\textrm{End}}
\newcommand{\Map}{\mathrm{Map}}
\newcommand{\spec}{\textrm{spec}}

\newcommand{\isom}{\cong}
\newcommand{\isoms}{\underline{\textrm{Isom}}}

\newcommand{\boper}{\mathcal{B}(\mathcal{H})}

\newcommand{\Ga}{\mathcal{G}}
\newcommand{\pr}{\mathrm{pr}}
\newcommand{\CAR}{\textrm{CAR}}
\newcommand{\I}{\mathcal{I}}
\newcommand{\Ab}{\mathfrak{Ab}}

\newcommand{\ad}{\mathrm{ad}\,}
\newcommand{\Ad}{\mathrm{Ad}\,}

\newcommand{\CCl}{\C\ell}
\newcommand{\coker}{\mathrm{coker}\,}
\newcommand{\vac}[1]{\vert #1\rangle}
\newcommand{\bra}[1]{\langle #1\vert}

\newcommand{\Sh}{\mathrm{Sh}}
\newcommand{\plim}{\varprojlim}

\newcommand{\Ext}{\mathrm{Ext}}
\newcommand{\detp}{\mathrm{det}_p}
\newcommand{\St}{\mathrm{St}}
\newcommand{\op}{\mathrm{op}}
\newcommand{\Gbun}{\Ga r}

\begin{document}

\title{Anomalies in gauge theory and gerbes over quotient stacks}%
\author{Vesa Tähtinen}%
\address{Department of Mathematics and Statistics, University of Helsinki,  P.O. Box 68 (Gustaf Hällströmin katu 2b),
FI-00014 Helsinki, Finland}%
\email{vesa.tahtinen@helsinki.fi}%

\thanks{The author would like to thank Finnish Academy of Science and Letters, Vilho, Yrjö and Kalle Väisälä Foundation for financial support}

\begin{abstract}
In Yang-Mills theory one is interested in lifting the action of the gauge transformation group $\Ga=\Ga(P)$ on the space of connection one-forms $\A=\A(P)$, where 
$P\To M$ is a principal $G$-bundle over a compact Riemannian spin manifold $M$,
to the total space of the
Fock bundle $\F\To\A$ in a consistent way with the second quantized Dirac operators $\hat{\Dirac_A},\,A\in\A$. In general, there is an obstruction to this
called the Faddeev-Mickelsson anomaly,
and to overcome this one has to introduce a Lie group extension $\hat{\Ga}$, not necessarily central, of $\Ga$ that acts in the Fock bundle. The Faddeev-Mickelsson anomaly
is then essentially the class of the Lie group extension $\hat{\Ga}$.

When $M=S^1$ and $P$ is the trivial $G$-bundle, we are dealing with $S^1$-central extensions of loop groups $LG$ as in \cite{PS}.
However, it was first noticed in the pioneering works of J.~Mickelsson, \cite{Mi} and L.~Faddeev, \cite{Fad} that when $\dim M>1$ the group multiplication in $\hat{\Ga}$ depends also on
the elements $A\in\A$ and hence is no longer an $S^1$-central extension of Lie groups.

We give a new interpretation of certain noncommutative versions of
Faddeev-Mickelsson anomaly (see for example \cite{Ra}, \cite{LMR} and \cite{AM}) and show that the analogous Lie group extensions $\hat{\Ga}$ can be replaced with a \emph{Lie groupoid} extension of the action Lie groupoid
$\A\rtimes\Ga$, where $\A$ is now some relevant abstract analog of the space of connection one-forms. Then at the level of Lie groupoids, this extension proves out to be an $S^1$-\emph{central} extension and hence one may apply the general theory of these
extensions developed by K. Behrend and P. Xu in \cite{BXu}. This makes it possible to consider the Faddeev-Mickelsson anomaly as the class of this Lie groupoid extension or equivalently
as the class of a certain differentiable $S^1$-gerbe over the quotient stack $[\A/\Ga]$. We also give examples from noncommutative gauge theory where our construction can be
applied. 

The construction may also be used to give a geometric interpretation of the (classical) Faddeev-Mickelsson anomaly in Yang-Mills theory when $\dim M=3$.

\end{abstract}

\maketitle

\tableofcontents

\section{Introduction}
\subsection{Obstruction to canonical quantization of fermions in Yang-Mills theory (a.k.a Faddeev-Mickelsson anomaly)}
\subsubsection{Dirac operators}
Suppose that $(M,g^M)$ is a compact oriented Riemannian spin manifold of dimension $d=2n+1$ without boundary and
let $S$ be the \emph{spin bundle} of the spin manifold $M$.

Let $G$ be a finite dimensional semi-simple compact Lie group and $\rho:G\To\Aut_\C(V)$ a unitary complex representation of $G$
with respect to an inner product $(\cdot,\cdot)_V$ on $V$, i.e.
$(\rho(g)x,\rho(g)y)=(x,y)$ for all $g\in G$ and $x,y\in V$. Next suppose that $\pi:P\To M$ is an arbitrary principal $G$ bundle and
form the associated vector bundle $E=P\times_\rho V$. One can show that since $\rho$ is unitary the associated vector bundle $E$ is
a Hermitean vector bundle with Hermitean metric $h^E$.

Denote by $\A$ the
space of $\Lie=\lie(G)$ valued connection $1$-forms on $P$ and by $\Ge$ the based gauge transformation group. It is known that
$\A/\Ge$ is a smooth infinite dimensional I.L.H. manifold, \cite{Pa}. To each $A\in\A$ one can associate a 
Dirac operator $\Dirac_A:\Gamma(\Ev)\To\Gamma(\Ev)$, where $\Ev:=S\tensor E$. This extends to an operator on
$\Hilb=\Ltwo(\Ev)$, the Hilbert space of square integrable sections of the vector bundle $\Ev$. The
domain of $\Dirac_A$ in $\Hilb$ is known to be $H^1(M;S)$, the first Sobolev space, \cite{Boss}.

One knows from functional analysis that $\Dirac_A$ is a
\emph{Fredholm} operator since it is elliptic and the manifold $M$ is compact. Thus $\dim\ker\Dirac_A<\infty$ and $\dim\coker\Dirac_A<\infty$. Moreover, the gauge transformation group $\Ge$ acts on $\Hilb$ and the
Dirac operator $\Dirac_A$ satisfies the following equivariance condition
$$
g\Dirac_A g^{-1}=\Dirac_{A^g}
$$
for all $g\in\Ge$.

\subsubsection{Fock bundle}

For each $A\in\A$ s.t. $0\notin\spec(\Dirac_A)$ the operator $\Dirac_A$ produces a decomposition
$$
\Hilb=\Hilb_+(A)\oplus\Hilb_-(A),
$$
where the spaces $\Hilb_\pm$ are the
corresponding eigenspaces to the
positive and negative eigenvalues of the Dirac operator $\Dirac_A$, respectively.
Corresponding to this decomposition
there exists an irreducible Dirac representation of the
representation of the algebra $\CAR(\Hilb)=:\CCl(\Hilb\oplus\bar{\Hilb})$ (the algebra of \emph{canonical anticommutation relations} or the
algebra of \emph{fermion fields}) on the
\emph{Fock space}
\begin{eqnarray}
\F_A&:=&\bigwedge\Big(\Hilb_+(A)\oplus\bar{\Hilb}_-(A)\Big)=\bigwedge\Hilb_+(A)\tensor\bigwedge\bar{\Hilb}_-(A)\nonumber\\
&=&\bigoplus_{p,q}\Big(\bigwedge^p\Hilb_+(A)\tensor\bigwedge^q\bar{\Hilb}_-(A)\Big),\nonumber
\end{eqnarray}
where physically the subspace $\bigwedge^p\Hilb_+(A)\tensor\bigwedge^q\bar{\Hilb}_-(A)$ consists of the states with $p$ particles and $q$ antiparticles, all of positive energy.
\footnote{Here $\bar{\Hilb}_-$ denotes the abstract complex conjugate space to $\Hilb_-$. It is a copy of $\Hilb_-$ with the scalars acting in a conjugate
way: $\lambda\cdot\bar{\xi}=(\lambda\cdot\xi)^{\!-}$; we don't suppose that there is a complex conjugation operation defined inside the Hilbert space $\Hilb$.}
A $\CAR$-representation $\psi_A:\CAR\To\End(\F_A)$ is  determined by giving a \emph{vacuum} vector $\vac{0_A}\in\F_A$
characterized by the property that
$$
\psi_A^\ast(u)\ket{0_A}=0=\psi_A(v)\ket{0_A},\quad\textrm{for all }u\in\Hilb_-(A),\,v\in\Hilb_+(A).
$$
\begin{defn}
Two representations of the $\CAR$-algebra are said to be
\emph{equivalent} if it is possible to represent them in the same Fock space in such a way that both corresponding vacuum vectors will be of finite norm.
\end{defn}
\begin{thm}
Two different polarizations $\Hilb=\Hilb_+\oplus\Hilb_-=W_+\oplus W_-$ define equivalent Dirac representations of the
$\CAR$-algebra if and only if the projections $\pr_{W_+}^-:W_+\To\Hilb_-$ and $\pr_{W_-}^+:W_-\To\Hilb_+$ are Hilbert-Schmidt.
\end{thm}
\begin{thm}[Shale-Stinespring]\label{ShaleStine}
Two Dirac representation of the $\CAR$-al\-ge\-bra defined by a pair of polarizations $\Hilb_+$ and $\Hilb_+'$
are equivalent if and and only if there is $g\in\U_{res}(\Hilb)$ such that $\Hilb_+'=g\cdot\Hilb_+$. In addition, in order that an element $g\in\U(\Hilb)$ is
\emph{implementable} in the Fock space, i.e. there is a unitary operator $\hat{g}\in\U(\F)$ such that
$$
\hat{g}\psi^\ast(v)\hat{g}^{-1}=\psi^\ast(gv),\quad\textrm{for all }v\in\Hilb,
$$
and similarly for the $\psi(v)$'s, one must have $g\in\U_{res}(\Hilb)$.
\end{thm}
Here
$\U_{res}(\Hilb)$ is the group of unitary operators $g$ in the polarized Hilbert space $\Hilb=\Hilb_+\oplus\Hilb_-$ such that the off-diagonal
blocks are Hilbert-Schmidt operators.

One would like to glue somehow the different $\CAR$-algebra representations $\F_A$ into an infinite-dimensional Hilbert bundle $\F$ over $\A$ with
a continuous section $s_\F:\A\To\F$ such that $s_\F(A)=\ket{0_A}$ (a Dirac
representation if fixed by a given vacuum vector so this way it is possible to define what
we mean by a continuously varying family of $\CAR$-representations). First, to construct a bundle of Fock spaces one can use the following trick:

One replaces the operator  $\Dirac_A$ with
the operator $\Dirac_A-\lambda$, where $\lambda\in\Real,\lambda\notin\spec(\Dirac_A)$. This way, one obtains a decomposition
$$
\Hilb=\Hilb_+(A,\lambda)\oplus\Hilb_-(A,\lambda),
$$
with the corresponding (irreducible) Fock space representation
$$
\rho_{A,\lambda}:\CAR(\Hilb)\To\End(\F_{A,\lambda})
$$
of the $\CAR$-algebra.

The Fock spaces $\F_{A,\lambda}$ depend on the choice of the \emph{vacuum level} $\lambda$. However, for $\lambda,\mu\notin\spec(\Dirac_A)$ there
exists a natural projective isomorphism
\begin{equation}\label{Fock}
\F_{A,\lambda}\equiv\F_{A,\mu}\mod\C^\times,
\end{equation}
allowing us to glue the different Fock spaces $\F_{A,\lambda}$ together into an infinite dimensional
\emph{projective} Fock bundle $\Proj\F$ over $\A$, \cite{Ar}. One can show that since $\A$ is contractible
as an affine space, there exists a trivial vector bundle $\F=\A\times\F_0$ over $\A$ whose projectivization is projectively isomorphic to $\Proj\F$.

Now the fibre of $\F$ at $A\in\A$ is equal to $\F_A\isom\F_0$ but
unfortunately for the energy polarization
$\Hilb=\Hilb_+(A)\oplus\Hilb_-(A)$
the map $A\mapsto\ket{0_A}$ does \emph{not} define a continuous section of $\F$ (or equivalently the map
$\A\To\Gr(\Hilb) :A\mapsto\Hilb_+(A)$ isn't continuous). This problem is resolved by intoducing another family $W(A)$ of polaritations
$\Hilb=W(A)\oplus W(A)^\bot$ parametrized by $A\in\A$ such that
\begin{enumerate}
\item The map $\A\To\Gr(\Hilb): A\mapsto W(A)$ \emph{is} continuous;
\item The corresponding CAR-algebra representations $\rho_A$  and $\rho_{W(A)}$ induced by the two polarizations are \emph{equivalent}.
\end{enumerate}
To construct such a family of polarizations one proceeds as follows (see \cite{Mi5} for details): Each $A\in\A$ defines a Grassmannian manifold $\Gbun_{res}(A)$ consisting
of all closed subspaces $W\subseteq\Hilb$ such that the difference $\pr_{\Hilb_+(A)}-\pr_W\in\mathcal{L}(\Hilb)$ is a Hilbert-Schmidt operator. One
can show that these spaces can be glued together to form a locally trivial fibre bundle over $\A$, called the \emph{Grasmannian} bundle $\Gbun$. The
question now is that does this bundle admit a global section $A\mapsto W(A)$? If it does the $W(A)$'s give us a family of polarizations with the required properties.

Luckily, the answer to our question is ``\emph{yes}''. This is because $\Gbun$ happens to be an associated bundle to an $\U_{res}(\Hilb)$-bundle $P\To\A$,
$$
\Gbun=P\times_{\U_{res}(\Hilb)}\Gr_{res}(\Hilb),
$$
where the fibre of $P$ at $A\in\A$ is
$$
P_A=\{g\in\U(\Hilb)\mid g\cdot\Hilb_+\in\Gbun_{res}(A)\}
$$
and
$\Gr_{res}(\Hilb)$ is the \emph{restricted} Grassmannian of Segal and Wilson (see Appendix A).
Now
$$
\Gr_{res}(\Hilb)\isom\U_{res}(\Hilb)/(\U(\Hilb_+)\times\U(\Hilb_-))
$$
and by a result of N.~Kuiper the subgroup $\U(\Hilb_+)\times\U(\Hilb_-)$ is contractible and so $\Gbun$ has
a global section if and only if $P$ is trivial. This happens to be the case since $\A$ is contractible as an affine space.
\subsubsection{Second quantizing gauge transformations}
After a certain necessary renormalization process, introduced by Mickelsson in \cite{Mi3}, on operations on the one-particle Hilbert space $\Hilb$ (e.g. the action of gauge transformation group)
one would hope to lift the action of $\Ga$ on $\A$ to an action on $\F$ so that the diagram
$$
\xymatrix{
\F\ar[r]^{\Gamma_A(g)}\ar[d] & \F\ar[d]\\
\A\ar[r]^{g} & \A
}
$$
commutes and
$$
\Gamma_A(g)\hat{\Dirac}_A\Gamma^{-1}_A(g)=\hat{\Dirac}_{A^g},
$$
where $\hat{\Dirac}_A$ is the second quantized Dirac operator.
Unfortunately, there is an obstruction to this. To study this, it is useful to switch to the Lie algebra picture.
\begin{defn}\label{secq}
Second quantization of an infinitesimal gauge transformation is the map $d\Gamma_A:\Dom(A)\subseteq\lie(\Ga)\To\End(\F_A)$ characterized by
\begin{eqnarray}
[d\Gamma_A(X),\psi_A^\ast(v)]&=&\psi_A^\ast(X\cdot v),\quad\textrm{for all } v\in\Hilb,\label{Bogo}\\
\bra{0_A}d\Gamma_A(X)\vac{0_A}&=&0.
\end{eqnarray}
\end{defn}
Here we may choose the domain $\Dom(A)$ of $d\Gamma_A(X)$ to be the set
$$
\Dom(A)=\{X\in\lie(\Ga)\mid[\epsilon_A,X]\textrm{ is Hilbert-Schmidt}\},
$$
where $\epsilon_A=\pm$ on $\Hilb_\pm(A)$. Moreover, supposing there exists a described lift $\Gamma_A:\Ga\To\End(\F)$ we should
have
$$
\Gamma_A(e^{iX})=e^{id\Gamma_A(X)},\quad\textrm{for all } X\in\lie(\Ga).
$$
In view of this, equation (\ref{Bogo}) can be written as
$$
\Gamma_A(e^{iX})\psi_A^\ast(v)\Gamma^{-1}_A(e^{iX})=\psi_A^\ast(e^{iX}\cdot v),\quad\textrm{for all } X\in\lie(\Ga),\,v\in\Hilb
$$
relating Definition \ref{secq} to Theorem \ref{ShaleStine}.

Next, we introduce the so called \emph{Gauss law generators} acting on (Schrödinger wave) functions $\phi:\A\To\Hilb$,
$$
G_A(X)=X+\mathcal{L}_X,
$$
where $A\in\A,\,X\in\lie(\Ga)$ and the \emph{Lie derivative} $\mathcal{L}_X$ is defined so that
$$
\Big(\mathcal{L}_X\phi\Big)(A)=\frac{d}{dt}\phi(A^{e^{tX}})\Big\vert_{\,t=0}
$$
Their second quantization is defined to be
$$
d\Gamma(G_A(X))=d\Gamma_A(X)+\mathcal{L}_X,
$$
where $X\in\lie(\Ga)$.
The renormalization procedure makes it possible to consider $d\Gamma_A(X)$ acting on $\F_0$ instead of $\F_A$. Now
the second quantized Gauss law generators do not have anymore the same Lie algebra bracket as $\lie(\Ga)$ but instead
$$
[d\Gamma(G_A(X)),d\Gamma(G_A(Y))]=d\Gamma([G_A(X),G_A(Y)])+c(X,Y;A),
$$
where $c(X,Y;A)$ is a $\Map(\A,\Real)$-valued Lie algebra cocycle of $\lie(\Ga)$ called the \emph{Schwinger term}. This is the sought obstruction term. The connection with 
\emph{bundle gerbes} comes from
a transgression map $\tau$,
$$
H^3(\A/\Ge,\Z)\To H^3(\A/\Ge,\Real)\isom H^3_{DR}(\A/\Ge)\buildrel\tau\over\To H^2(\lie(\Ga),\Map(\A,\Real))
$$
studied in \cite{CMW}.

In \cite{CMM} Carey, Mickelsson and Murray constructed explicitly the bundle gerbe in question using a collection of local determinant line bundles on the 
smooth Fréchet manifold $\A/\Ge$ that satisfy certain compatibility conditions.
Let us recall this construction briefly.

Define for all $\lambda\in\Real$ the open subsets
$$
U_\lambda=\{A\in\A\mid\lambda\notin\spec(\Dirac_A)\}\subseteq\A.
$$
These form an open cover for $\A$.
Over each intersection $U_{\lambda\mu}:=U_\lambda\cap U_\mu$ there exists a line bundle $\Det_{\lambda\nu}$, whose fibre
$\Det_{\lambda\nu}(A)$ at $A\in\A$ is related to (\ref{Fock}) by the equation
$$
\F_{A,\lambda}=\Det_{\lambda\mu}(A)\tensor\F_{A,\mu}
$$
(thus giving the phase) and defined so that
$$
\Det_{\lambda\mu}(A)=\bigwedge^{max}\big(\Hilb_+(A,\lambda)\cap \Hilb_-(A,\mu)\big)
$$
for $\lambda<\mu$ and $\Det_{\mu\lambda}:=\Det_{\lambda\mu}^{-1}$.
The phase is related to the arbitrariness in filling the Dirac sea between
vacuum levels $\lambda$ and $\mu$. Such a filling corresponds to an exterior product
$v_1\wedge v_2\wedge\ldots\wedge v_m$ of a complete orthonormal set of eigenvectors
$\Dirac_A v_i=\lambda_i v_i$ with $\lambda<\lambda_i<\mu$. A rotation of the eigenvector basis gives a multiliplication
of the exterior product by the determinant of the rotation.
Now, since the exterior product satisfies the 'exponential law'
$$
\bigwedge^{max}(V\oplus W)=\bigwedge^{max} V\tensor\bigwedge^{max} W
$$
for finite dimensional vector spaces $V$ and $W$,
one sees that over the triple intersections
$U_{\lambda\lambda'\lambda''}:=U_\lambda\cap U_{\lambda'}\cap U_{\lambda''}$
$$
\Det_{\lambda\lambda'}\tensor\Det_{\lambda'\lambda''}=\Det_{\lambda\lambda''},
$$
so that the collection $\{\Det_{\lambda\mu}\}$ of local line bundles define a \emph{bundle gerbe} on $\A$. These local determinant line bundles are actually $\hat{\mathcal{G}}$-equivariant,
where $\hat{\Ga}$ is the group extension of $\Ga$ integrating the Lie algebra extension of $\lie(\Ga)$ determined by the Scwhinger term, and so
descend to the moduli space $\A/\Ge$ giving us the bundle gerbe whose Dixmier-Douady class transgresses to the Schwinger term.

¨\subsection{Main results} We use \emph{differentiable gerbes} of Behrend and Xu  \cite{BXu} instead of bundle gerbes to describe 
geometrically the noncommutative version of Faddeev-Mickelsson anomaly. This allows us to consider situations
where a relevant generalized gauge transformation group $\Ga$ (e.g. $\U_p(\Hilb)$) no longer acts 
freely and transitively on some space of generalized connection one-forms $\A$ (e.g. $\Gr_p(\Hilb)$). This is often the
case with noncommutative gauge theories, where it is hard to find a relevant gauge group acting nicely enough. 

In this picture the noncommutative
Faddeev-Mickelsson anomaly is given by the gerbe class $\omega\in H^2([\A/\Ga],\underline{S}^1)$ of a certain $S^1$-gerbe over the quotient stack $[\A/\Ga]$ or equivalently by the class of
a certain $S^1$-Lie groupoid extension of the action groupoid $\A\rtimes\Ga$ which we construct. When $\A/\Ga$ exists as a nice manifold (e.g. a Banach or an I.L.H. manifold) 
satisfying the smooth partition of unity property one knows that $[\A/\Ga]\isom\A/\Ga$ and
$H^2([\A/\Ga],\underline{S}^1)\isom H^2(\A/\Ga,\underline{S}^1)\isom H^3(\A/\Ga,\Z)$, where the last cohomology group classifies bundle gerbes, \cite{St}.

It was proven in \cite{LaMi} that in dimension equal to three and at the level of Lie group entensions one can revive the actual Faddeev-Mickelsson anomaly in (classical) Yang-Mills theory from a noncommutative
Faddeev-Mickelsson anomaly. Namely, one can pull-back the noncommutative Faddeev-Mickelsson anomaly Lie algebra cocycle and it proves out that this
represents the same class as the original Faddeev-Mickelsson anomaly cocycle. Hence our methods may also be used to describe the original Faddeev-Mickelsson anomaly on
a compact Riemannian spin manifold $M$, when $\dim M=3$.

\subsection*{Acknowledgments} The author would like to thank Professor Jouko Mickelsson for introducing the problem and giving many helpful comments. The work
was financially supported by the Finnish Academy of Science and Letters, Vilho, Yrjö and Kalle Väisälä Foundation. The author would also like to thank
Erwin Schrödinger International Institute for Mathematical Physics for hospitality where the work was initiated in summer 2006.

\section{NCG field theory examples}
Here we give two examples from noncommutative gauge theory in which it is difficult to find any relevant gauge transformation group $\Ga$ acting freely and transitively on the
space of connections $\A$. In that case $[\A/\Ga]$ is no longer a smooth manifold but rather a (differentiable) stack and hence bundle gerbes on it are not defined anymore. 
However, the author thinks one
might be able to develop some $\Ga$-equivariant bundle gerbe approach to Faddeev-Mickelsson anomalies in this setting, but since we actually work at the level of
Lie groupoids we prefer to speak about quotient stacks instead in the spirit of \cite{BXu}.
\subsection{Universal Yang-Mills theory of Rajeev}
Here we follow \cite{MiRa}, \cite{Ra} and \cite{Mi1}.
\subsubsection{Generalized Fredholm determinants}
Let $\Hilb$  be a complex infinite dimensional separable Hilbert space with
a given polarization $\Hilb=\Hilb_+\oplus\Hilb_-$. Let $\Lp^p$, where $p\geq 1$, denote the \emph{Schatten ideal}, i.e.
the space of linear operators $A:\Hilb\To\Hilb$ s.t.
$$
\norm{A}^p_p=\tr(A^\ast A)^{p/2}<\infty.
$$
Each $\Lp^p$ is a complete metric space with respect to the norm $\norm{\cdot}_p$.

Now for each $A\in\Lp^p$ define
$$
R_p(A)=-1+(1+A)\exp\Big[\sum_{j=1}^{p-1}(-1)^j\frac{A^j}{j}\Big].
$$
\begin{defn}[Generalized Fredholm determinants]
Let $A\in\Lp^p$ and define
$$
\detp(1+A):=\det(1+R_p(A)).
$$
\end{defn}
We have the following formula
$$
\log\detp(1+A)=\tr\Big((-1)^p\frac{A^p}{p}+(-1)^{p+1}\frac{A^{p+1}}{p+1}+\cdots\Big)
$$
so that $\log\detp(1+A)$ can be thought of as a \emph{regularization} of $\det(1+A)$, where the first
$p-1$ terms have been subtracted in the expansion of $\log(1+A)$.

The regularized determinants are \emph{not} multiplicative but instead we have the following proposition
\begin{prop}
For each $p\in\N^+$ there is a symmetric polynomial $\gamma_p(A,B)$ of two variables
$A,B\in 1+\Lp^p$ such that
$$
\detp\, AB=\detp A\cdot\detp B\cdot e^{\gamma_p(A,B)}.
$$
\end{prop}
\begin{defn}
$\omega_p(A,B)=\detp B\cdot e^{\gamma_p(A,B)}$.
\end{defn}

When $A$ is invertible it is known that
$$
\omega_p(A,B)=\frac{\det_p AB}{\det_p A}.
$$
More over, for $A,B,C\in1+\Lp^p$
$$
\omega_p(A,BC)=\omega_p(AB,C)\cdot\omega_p(A,B).
$$
\subsubsection{Generalized determinant line bundles}
Let $\Hilb$  be a complex infinite dimensional separable Hilbert space with
a given polarization $\Hilb=\Hilb_+\oplus\Hilb_-$. We
fix an orthonormal basis $\{e_n\}_{n\in\Z}$ of $\Hilb$ such that $e_n\in\Hilb_+$ for $n>0$ and $e_n\in\Hilb_-$ for $n\leq 0$.

Let $\GL_p(\Hilb)$ denote the group consisting of all invertible bounded linear operators of the form
$$
\Mat{a}{b}{c}{d},
$$
where $a:\Hilb_+\To\Hilb_+,\,d:\Hilb_-\To\Hilb_-,\,c:\Hilb_+\To\Hilb_-$ and $b:\Hilb_-\To\Hilb_+$ are linear
operators such that $b,c\in\Lp^{2p}$. The group $\GL_p(\Hilb)$ has a natural metric topology defined by
$$
d(g,g')=\norm{a-a'}+\norm{d-d'}+\norm{b-b'}_{2p}+\norm{c-c'}_{2p}.
$$
This makes $\GL_p(\Hilb)$ into a Banach-Lie group.
\begin{defn}[Grassmannian manifold]
Let $B_p(\Hilb)$ be the (closed) normal subgroup of the block triangular operators in $\GL_p(\Hilb)$ with $c=0$. Define the infinite-dimensional $p\,$:th
\emph{Schatten Grassmannian} by
$$
\Gr_p(\Hilb):=\GL_p(\Hilb)/B_p(\Hilb).
$$
As a homogeneous space of a Banach-Lie group, $\Gr_p(\Hilb)$ is a Banach-Lie group.
\end{defn}

The points of $\Gr_p$ can be thought of as infinite-dimensional closed subspaces $W\subseteq\Hilb$ such that
\begin{enumerate}
\item The projection $\pr_{\Hilb_+}:W\To\Hilb_+$ is a Fredholm operator;
\item The projection $\pr_{\Hilb_-}:W\To\Hilb_-$ belongs to the Schatten ideal $\Lp^{2p}$.
\end{enumerate}

\begin{defn}
A basis $w=\{w_n\}_{n=1,2,\ldots}$ of $W\in\Gr_p$ is said to be admissible (with respect to the basis $\{e_n\}_{n>0}$ of $\Hilb_+$)
if $w_+-1\in\Lp^p$, where $w_+$ is the (infinite) matrix defined by
$$
\pr_{\Hilb_+}w_i=\sum_{j>0}(w_+)_{ji}e_j.
$$
\end{defn}

\begin{defn}
Let
$$
\Ev_p:=\{(g,q)\mid g\in\GL_p,\, q\in\GL(\Hilb_+),\,aq^{-1}-1\in\Lp^p\}\subseteq\GL_p\times\GL(\Hilb_+),
$$
where $g=\Mat{a}{b}{c}{d}$, be the group whose group multiplication is given by
$$
(g_1,q_1)(g_2,q_2)=(g_1g_2,q_1q_2)
$$
and topology
by the norm
$$
\norm{(g,q)}=\norm{a}+\norm{d}+\norm{b}_{2p}+\norm{c}_{2p}+\norm{a-q}_p.
$$
Then $\Ev_p$ is a Banach-Lie group.
\end{defn}

\begin{defn}
Define $\GL^p=\GL(\Hilb_+)\cap(1+\Lp^p)$, where $p\in\N\cup\{\infty\}$;
$\Lp^0=\{\textrm{finite rank operators}\},\quad \Lp^\infty=\{\textrm{compact operators}\}$.
\end{defn}

\begin{defn}[Stiefel manifolds]
The infinite-dimensional $p\,$:th \emph{Schatten-Stiefel manifold}
$$
\St_p:=\Ev_p/B_p,
$$
where the action of $k=\Mat{\alpha}{\beta}{0}{\gamma}\in B_p$ is given by
$$
(g,q)\cdot k=(gk,q\alpha).
$$
\end{defn}

The Stiefel manifold $\St_p$ parametrizes all admissible basis of all infinite-di\-men\-sio\-nal planes $W\in\Gr_p$, see \cite{Mi1}.
It is in a natural way a principal $\GL^p$-bundle over $\Gr_p$, the $\GL^p$ action being given by
the basis transformations and the canonical projection $\St_p\To\Gr_p$ is chosen to be the mapping
associating to the basis $w$ the plane $W$ spanned by the vectors in $w$.

\begin{defn}[Generalized determinant line bundles]
Let
$$
\Det_{\!p}:=(\St_{\!p}\times\C)/\GL^p,
$$
where the right action of $\GL^p$ on $\St_p\times\C$ is defined so that
$$
(w,\lambda)\cdot t=(wt,\lambda\omega_p(w_+,t)^{-1}).
$$
\end{defn}

One can show that $\Det_p$ is a holomorphic line bundle over $\Gr_p$ where the projection map
is given by $[(w,\lambda)]\mapsto$ the plane spanned by $\{w_1,w_2,\ldots\}$. Moreover,
the group $\GL_p$ acts on the base manifold $\Gr_p$ but the action doesn't lift to the bundle
$\Det_p$ for $p\geq 1$.

Naturally there is also the dual determinant line bundle $\Det_p^\ast\To\Gr_p$.

\begin{lem}
Sections of $\Det_p^\ast$ can be identified with functions $\psi:\St_p\To\C$ such that
$$
\psi(wt)=\psi(w)\omega_p(w,t),\quad t\in\GL^p.
$$
\end{lem}

\subsubsection{The Abelian extension of $\GL_p$}

\begin{lem}
There are smooth functions $\alpha(g,q;w)$ on $\Ev_p\times\St_p$ s.t.
$$
\frac{\alpha(g,q;wt)}{\alpha(g,q;w)}=-\frac{\omega_p(w_+,t)}{\omega_p((gwq^{-1})_+,qtq^{-1})}.
$$
\end{lem}

\begin{thm}[Mickelsson and Rajeev, \cite{MiRa}]
Let $\Hilb$  be a complex infinite dimensional separable Hilbert space with
a given polarization $\Hilb=\Hilb_+\oplus\Hilb_-$. There is an Abelian extension of $\GL_p=:\GL_p(\Hilb)$
by $\Map(\Gr_p,\C^\ast)$ which acts on $\Det_p$. The extension is
$$
\widehat{\GL_p}=(\Ev_p\times\Map(\Gr_p,\C^\ast))/N,
$$
where $N$ is the normal subgroup consisting of elements $(1,q,\mu_q)$, where
$\mu_q(w)=\alpha(1,q,w)^{-1}\cdot\omega_p(w_+,q^{-1})^{-1},\,q\in\GL^p$.
\end{thm}

\begin{rem}
As a corollary, one obtains the Abelian Lie group extension $\widehat{\U}_p(\Hilb)$ of $\U_p(\Hilb)\subseteq\GL_p(\Hilb)$
by the group $\Map(\Gr_p,\C^\ast)$ by restriction.
\end{rem}

\subsubsection{Canonical formalism for universal gauge theory}
The configuration space in \emph{Universal Yang-Mills theory} is by definition
$$
\tilde{\A}=\Big\{\textnormal{bounded Hermitean } \tilde{A}:
\Hilb\To\Hilb\, \Big\vert\, \tilde{A}\in\Mat{\Lp^p}{\Lp^{2p}}{\Lp^{2p}}{\Lp^p}\Big\}.
$$
The subgroup $\U_p\subseteq\GL_p$ of unitaries plays the role of the gauge transformation group acting on the manifold $\tilde{\A}$ by the rule
$$
\tilde{A}\mapsto\tilde{g}\tilde{A}\tilde{g}^{-1}+\tilde{g}[\epsilon,\tilde{g}^{-1}].
$$
The operator $\tilde{g}[\epsilon,\tilde{g}^{-1}]$ is indeed of type
$$
\Mat{\Lp^p}{\Lp^{2p}}{\Lp^{2p}}{\Lp^p}
$$
since we know that for Schatten ideals
$$
\Lp^p\cdot\Lp^q\subseteq\Lp^r,
$$
where
$
1/r=1/p+1/q.
$

The space of ``electric fields'' is
$$
\tilde{\mathcal{E}}=\Big\{\textnormal{bounded Hermitean } \tilde{E}:
\Hilb\To\Hilb\, \Big\vert\, \tilde{E}\in\Mat{\Lp^{p/(p-1)}}{\Lp^{2p/(2p-1)}}{\Lp^{2p/(2p-1)}}{\Lp^{p/(p-1)}}\Big\}.
$$

The phase space of universal Yang-Mills theory is defined to be the direct sum $\tilde{\A}\oplus\tilde{\mathcal{E}}$.
This space has a natural exterior derivative operator $\tilde{d}:\tilde{\A}\oplus\tilde{\mathcal{E}}\To\tilde{\A}\oplus\tilde{\mathcal{E}}$,
$$
\tilde{d}(\tilde{A},\tilde{E}):=([\epsilon,\tilde{E}],[\epsilon,\tilde{A}]_+),
$$
where $[\cdot,\cdot]_+$ means the anti-commutator. The elements of the form $(\tilde{A},0)\in\tilde{\A}\oplus\tilde{\mathcal{E}}$
are said to be of \emph{odd} degree and respectively the elements of the form $(0,\tilde{E})\in\tilde{\A}\oplus\tilde{\mathcal{E}}$
are said to be of \emph{even} degree. Clearly, $\tilde{d}$ maps even operators to odd operators and vice versa. Furthermore,
$\tilde{d}^2(\tilde{A},\tilde{E})=0$, since $\epsilon^2=1$.

The exterior derivative operator $\tilde{d}$ makes it possible to define the \emph{curvature} $\tilde{F}$ for every $\tilde{A}\in\tilde{\A}$,
$$
\tilde{F}:=\tilde{d}\tilde{A}+\tilde{A}^2.
$$
This is
an even operator in the sense we just defined. The curvature transforms covariantly
under gauge transformation, $\tilde{F}\mapsto\tilde{g}\tilde{F}\tilde{g}^{-1}$.
\begin{defn}
We say that a generalized vector potential/connection $1$-form $\tilde{A}\in\tilde{\A}$ is \emph{flat} if its curvature
$\tilde{F}=0$.
\end{defn}
\begin{prop}
The space of flat connections in universal Yang-Mills theory with gauge transformation group $\U_p(\Hilb)$ can be identified with
the $p\,$:th Schatten Grassmannian
$$
\Gr_p(\Hilb)\isom\U_p/(\U(\Hilb_+)\times\U(\Hilb_-)).
$$
\end{prop}
\subsubsection{Generalized Fock bundles over $\Gr_2(\Hilb)$}
An excellent reference for this subsection is \cite{Mi4}.

First, recall from \cite{PS} the geometric construction of the Fermionic Fock space
as the space of holomorphic sections of a complex line bundle $\Det_1^\ast$ over $\Gr_1$.
We want to generalize this to higher dimensional cases.

We suppose our Schatten Grassmannian $\Gr_2(\Hilb)$ is defined by a splitting $\Hilb=\Hilb_+\oplus\Hilb_-$.
\begin{defn}
Let $F\in\Gr_2(\Hilb_+\oplus\Hilb_-)$
and let $\Hilb=F\oplus F^\bot$ be the associated splitting. We define
the \emph{generalized} Fock space $\F_{F}$ by
$$
\F_{F}:=\Gamma(\Det_2^\ast(F\oplus F^\bot)),
$$
where $\Det_2^\ast(F\oplus F^\bot)\To\Gr_2(F\oplus F^\bot)$ is the dual of the $2\,$:nd determinant line bundle $\Det_2(F\oplus F^\bot)$.
\end{defn}

Now the problem with the above construction is that the Fock spaces $\F_{F}$ depend on a choice of admissible basis $f$ in each $F\in\Gr_2(\Hilb_+\oplus\Hilb_-)$:
\begin{lem}
Fix an admissible basis $f=\{f_1,f_2,\ldots\}$ of $F\in\Gr_2(\Hilb_+\oplus\Hilb_-)$.
Then a section $\tilde{\psi}_F\in\Gamma(\Det_2^\ast(F\oplus F^\bot))$ can be identified with a function $\psi_F:\St_{\,2}(F\oplus F^\bot)\To\C$ satisfying
\begin{equation}\label{gfock}
\psi_F(wt)=\psi_F(w)\cdot\omega_2(w^{(f)},t),\quad t\in\GL^2(F\oplus F^\bot),
\end{equation}
where $w(f)$ is the matrix relating the $F$-projection to the basis $\{f_n\}$, i.e.
$$
\pr_F(w_n)=\sum_j w_{jn}^{(f)} f_j
$$
and
$$
\omega_2(w^{(f)},t)=\frac{\det_2 w^{(f)}t}{\det_2 w^{(f)}}.
$$
\end{lem}

In fact, what we have consructed is a fibre bundle over $\St_{\,2}(\Hilb_+\oplus\Hilb_-)$ and \emph{not} over $\Gr_2(\Hilb_+\oplus\Hilb_-)$. We need to modify the
situation a bit to obtain a bundle over $\Gr_2(\Hilb_+\oplus\Hilb_-)$ and for this we proceed as follows.

Since the definition of a section $\psi$ depends on $f$ we shall write explicitly $\psi=\psi(w,f)$ and consider
these also as functions of $f$.

\begin{prop}
Functions $\psi_F:\St_{\,2}(F\oplus F^\bot)\times\St_{\,2}(F\oplus F^\bot)\To\C$ satisfying equation (\ref{gfock}) and
\begin{equation}
\psi_F(w,ft)=\psi_F(w,f)\cdot\omega_2(w^{(f)},t^{-1}),\quad t\in\GL^2(F\oplus F^\bot)
\end{equation}
can be identified with sections of a vector bundle $\F'$ over $\Gr_2(\Hilb_+\oplus\Hilb_-)$ which is a tensor product of the
determinant bundle $\Det_2(\Hilb_+\oplus\Hilb_-)$ and a trivial Fock bundle $\mathcal{B}$ (with fibre $\F_{\Hilb_+}$) over $\Gr_2(\Hilb_+\oplus\Hilb_-)$.
\end{prop}

\begin{defn}
We define the \emph{generalized} Fock bundle $\F'$ over $\Gr_2(\Hilb_+\oplus\Hilb_-)$ by
$$
\F':=\mathcal{B}\tensor\Det_2.
$$
\end{defn}

Motivated by this, one may define the obstruction to canonical quantization in universal Yang-Mills theory to be the
class of the Abelian Lie group extension $\widehat{\U}_2(\Hilb)\To\U_2(\Hilb)$. 
\subsection{NCG theory model of Langmann, Mickelsson and Rydh}
Our references in this section are \cite{LMR}, \cite{VG} and \cite{Co}.
\subsubsection{The space of generalized vector potentials}
Let $(\Hilb,D_0)$ be a tame $p^+$ summable $K$-cycle over the $^\ast$-algebra
$$
\Al=\{A\in\boper\mid
[\abs{D_0},A]\in\Lp^{p+}(\Hilb),\,[D_0,A]\in\boper\}
$$
with
$\pi:\Al\To U(\Hilb)$ the corresponding unitary representation and
$\Gamma:\Hilb\To\Hilb$ a grading operator. Denote by $\epsilon=D_0/\abs{D_0}$ the sign of the (abstract) Dirac operator.
Using the representation $\pi$, the equivalence classes $\alpha\in\A:=\Omega_{D_0}^1(\Al)$ can then be presented in the form
$$
\alpha=a_0[D_0,a_1],\quad a_0,a_1\in\Al,
\quad
\textrm{or}
\quad
\alpha=a_0[\epsilon,a_1],\quad a_0,a_1\in\Al.
$$
It follows that all the operators $\alpha\in\A$ satisfy the condition $[\epsilon,\alpha]\in\boper$.
\subsubsection{Gauge transformation group}
We assume our
(Hermitean) vector bundle $\Ev$ on $\Al$ to be trivial and of rank one, i.e. $\Ev=\Al$. Hence the gauge group $\U(\Ev)$ is given by
$$
\U(\Ev)=\U_{p+}=\{u\in\Al\mid uu^\ast=u^\ast u=1\}.
$$
Any element $g\in \U_{p+}(\Hilb)$ satisfies $[\epsilon,g]\in\Lp^{p+}$. This is seen to be the group of unitaries in the group
$$
\GL_{p+}:=\{g\in\Al\mid g \textrm{ is invertible}\}.
$$

\subsubsection{Family of (abstract) Dirac operators over $\A$}
We consider bounded perturbations $D_A$ of the `\emph{free Dirac operator}' $D_0$ that are of the form
$D_A=D_0+A$, where $A\in\A$ and the sign operator $F_A:=D_A/\abs{D_A}$ satisfies
$$
F_A=F_A^\ast=F_A^{-1}\in\boper,\quad
F_A-\epsilon\in\Lp^{p+}.
$$

Following the ideas of \cite{La1} and \cite{La2}, one can see that the sign operator $F_A$ can thus be thought as an element of
the \emph{weak-}$\Lp^p$ Grassmannian $\Gr_{p+}(\Hilb)$  defined analogously with the Schatten Grassmannian $\Gr_p(\Hilb)$ except
that now we require that the projection $\pr_{\Hilb_-}:W\To\Hilb_-$ belongs to the weak-$\Lp^p$ space $\Lp^{p+}$
instead of the Schatten ideal $\Lp^p$. More over, the Grassmannian $\Gr_{p+}$ has a natural action of the group $\GL_{p+}$.

This motivates us to consider the obstruction of canonically quantizing fermions in this NCG gauge theory model as the class of the
group extension $\hat{\U}_{p+}$  acting on the total space of the determinant line bundle
$\Det_{p+}\To\Gr_{p+}$ analogously with what we did in the case of universal Yang-Mills theory.

The group extension $\widehat{\GL}_{p+}$ can be constructed in the same vein as in \cite{AM}. However, one
has to pay attention to the properties of generalized traces, \cite{LMR}.
\section{Differentiable $S^1$-gerbes and $S^1$-Lie groupoid central extensions}
The main reference in this section is \cite{BXu}.
\subsection{Stacks}
Let $\Sc$ be either the category of all finite dimensional $\Cinfty$-manifolds with $\Cinfty$-maps as morphisms,
or the category of all (infinite dimensional) $\Cinfty$-Banach manifolds with the corresponding smooth maps. We endow $\Sc$ with the Grothendieck topology, whose
covering families $\{U_i\To X\}$ are
local diffeomorphisms
$U_i\To X$
such that the total map
$\coprod_i U_i\To X$ is surjective.

\begin{defn}
A \emph{category fibered in groupoids} $\stack\To\Sc$ is a category $\stack$, together with a functor
$\pi:\stack\To\Sc$, such that the following two conditions are satisfied:
\begin{enumerate}
\item For every arrow $V\To U$ in $\Sc$, and every object $x$ of $X$ lying over $U$, $\pi(x)=U$, there exists an arrow $y\To x$ in $\stack$ lying
over $V\To U$, i.e. $\pi(y\To x)=V\To U$.
\item For every commutative diagram $W\To V\To U$ in $\Sc$ and arrows $z\To x$ lying over $W\To U$ and $y\To x$, there exists a unique
arrow $z\To y$ lying over $W\To V$, such that the composition $z\To y\To x$ equals $z\To x$.
\end{enumerate}
\end{defn}

\begin{ex}
Manifolds $X\in\Ob(\Sc)$ give groupoid fibrations. To see this, let $\underline{X}$ denote the category where
$$
\Ob(\underline{X})=\{(S,u)\mid S\in\Ob(\Sc),u\in\Hom_\Sc(S,X)\}
$$
and a morphism $(S,u)\To(T,v)$ of objects is a morphism $f:S\To T$ such that $u=v\circ f$, i.e. an $X$ morphism.
\end{ex}

\begin{defn}
Let $\pi:\stack\To\Sc$ be a category fibered in groupoids. Then $\stack$ is called a \emph{stack} over $\Sc$ if the following three axioms are
satisfied:
\begin{enumerate}
\item For any $\Cinfty$ manifold $X\in\Ob(\Sc)$, any two objects $x,y\in\Ob(\stack)$ lying over $X$, and any two isomorphims
$\phi,\psi:x\To y$ over $X$ such that $\phi\vert U_i=\psi\vert U_j$ for all $U_i$ in a covering family $\{U_i\To X\}$, then $\phi=\psi$.
\item For any $X\in\Ob(\Sc)$, any two objects $x,y\in\Ob(\stack)$ lying over $X$, a covering family $\{U_i\To X\}$, and a collection of
isomorphisms $\phi_i: x\vert U_i\To y\vert U_i$ such that $\phi_i\vert U_i\times_X U_j=\phi_j\vert U_i\times_X U_j$ for all $i,j$, there exists
an isomorphism $\phi: x\To y$ such that $\phi\vert U_i=\phi_i$ for all $i$.
\item For every $X\in\Ob(\Sc)$, every covering family $\{U_i\To X\}$, every family $\{x_i\}$ of objects $x_i$ in the fibre
$\stack_{U_i}$, and every family of morphims $\{\phi_{ij}\}$, $\phi_{ij}:x_i\vert U_i\times_X U_j\To x_j\vert U_i\times_X U_j$ satisfying the cocycle condition
$\phi_{jk}\circ\phi_{ij}=\phi_{ik}$ in the fibre $\stack_{U_i\times_X U_j\times_X U_k}$, there exists an object $x$ over $X$, together with isomorphisms
$\phi_i:x\vert U_i\To x_i$ such that $\phi_{ij}\circ\phi_i=\phi_j$ over $U_{ij}$.
\end{enumerate}
\end{defn}
\begin{rem} Here condition (2) means that morphisms glue and condition (3) says that objects glue (descent data is effective).
Conditions (1) and (2) imply that
for fixed $X\in\Ob(\Sc)$, $x,y\in\stack_X$,
$\isoms(x,y)$ is a sheaf on $\Sc/X$.
\end{rem}

The morphisms of stacks are morphisms of their underlying groupoid fibrations.

\begin{ex}[Manifolds]
For every manifold $X\in\Ob(\Sc)$ the groupoid fibration $\underline{X}$ is a stack.
\end{ex}

\begin{ex}[Quotient stacks]
Let $G\in\Ob(\Sc)$ be a Lie group acting on a manifold $X\in\Ob(\Sc)$. Define the \emph{quotient stack} $[X/G]$ as the
category whose objects are principal $G$-bundles $\pi:P\To S$, where all manifolds and structure maps are in $\Sc$, together with a $G$-equivariant morphism
$\alpha\in\Hom_\Sc(P,X)$. A morphism in $[X/G]$ is a Cartesian diagram in $\Sc$
$$
\xymatrix{
P'\ar[r]^p \ar[d]_{\pi'} & P\ar[d]^\pi\\
S'\ar[r]^f & S
}
$$
such that $\alpha\circ p=\alpha'$. The projection functor $\pi_{[X/G]}:[X/G]\To\Sc$ associates to a principal
$G$-bundle $\pi:P\To S$ its base space $S$ and to a morphism as above the map $f:S'\To S$ in $\Sc$.
Choosing $X=\bullet$, a point, one obtains the \emph{classifying stack} $BG$.

If $G$ acts properly and freely, i.e. $X \To X/G$ is a $G$-bundle,
then $[X/G]\isom X/G$, see \cite{He}, Remark 1.6.
\end{ex}


\begin{defn}
A stack $\stack$ over $\Sc$ is called \emph{differentiable} or a $\Cinfty$ stack, if there exists a manifold
$X\in\Ob(\Sc)$ and a surjective representable submersion $x:X\To\stack$. In this case $X$ together with the structure morphism $x$
is called an \emph{atlas} for $\stack$ or a \emph{presentation} of $\stack$.
\end{defn}
\begin{ex}[Quotient stacks] An atlas is given by the quotient map $X\To [X/G]$, defined by
the trivial $G$-bundle $G\times X\To X$ and $\alpha:G\times X\To X$ being the action map.
\end{ex}

\subsection{Lie groupoids}

\begin{defn}
A \emph{Lie groupoid} $\Gamma=X_1\rightrightarrows X_0$ consists of
\begin{itemize}
\item Two smooth manifolds $X_1\in\Ob(\Sc)$ (the \emph{morphisms} or \emph{arrows}) and $X_0\in\Ob(\Sc)$ (the \emph{objects} or \emph{points});
\item Two smooth surjective submersions $s:X_1\To X_0$ the \emph{source} map and $t:X_1\To X_0$ the \emph{target} map;
\item A smooth embedding $e:X_0\To X_1$ (the \emph{identities} or \emph{constant arrows});
\item A smooth involution $i:X_1\To X_1$, (the \emph{inversion}) also denoted $x\mapsto x^{-1}$;
\item A multiplication
$$
m:\Gamma^{(2)}\To\Gamma,
$$
$$
(x,y)\mapsto x\cdot y,
$$
where $\Gamma^{(2)}=X_1\times_{s,t} X_1=\{(x,y)\in X_1\times X_1\mid s(x)=t(y)\}$. Notice, that $\Gamma^{(2)}$ is a smooth manifold, since
$s$ and $t$ are submersions. We require the multiplication map $m$ to be smooth and that
\begin{enumerate}
\item $s(x\cdot y)=s(y),\quad t(x\cdot y)=t(x)$,
\item $x\cdot(y\cdot z)=(x\cdot y)\cdot z$,
\item $e$ is a section of both $s$ and $t$,
\item $e(t(x))\cdot x=x=x\cdot e(s(x))$,
\item $s(x^{-1})=t(x),\quad t(x^{-1})=s(x)$,
\item $x\cdot x^{-1}=e(t(x)),\quad x^{-1}\cdot x=e(s(x))$,
\end{enumerate}
whenever $(x,y)$ and $(y,z)$ are in $\Gamma^{(2)}$.
\end{itemize}
\end{defn}

\begin{rem}
When $\Sc$ is the category of smooth Banach manifolds, we call $\Gamma=X_1\rightrightarrows X_0$ a Banach-Lie groupoid.
\end{rem}
\begin{defn}
A morphism of Lie groupoids $(\Psi,\psi):[X_1'\rightrightarrows X_0']\To[X_1\rightrightarrows X_0]$
are the following commutative diagrams:
$$
\xymatrix{
X_1'\ar@<2pt>[d]^{t'}\ar@<-2pt>[d]_{s'}\ar[r]^\Psi & X_1\ar@<2pt>[d]^{t}\ar@<-2pt>[d]_{s}& &X_1\ar[r]^\Psi & X_1'\\
X_0'\ar[r]^{\psi} & X_0 & & X_0'\ar[u]^{e'}\ar[r]^\psi & X_0\ar[u]_e
}
$$
$$
\xymatrix{
X_1'\times_{s',t'} X_1'\ar[r]^{\Psi\times\Psi}\ar[d]_{m'}& X_1\times_{s,t}X_1\ar[d]^m & & X_1'\ar[r]^\Psi\ar[d]_{i'} & X_1\ar[d]^{i}\\
X_1'\ar[r]^\Psi & X_1 & & X_1'\ar[r]^\Psi & X_1
}
$$
\end{defn}

\begin{ex}
A Lie group $G$ is a Lie groupoid over a point, $G\rightrightarrows\bullet$.
\end{ex}

\begin{ex}
Let $M$ be a differentiable manifold and $G$ a Lie group acting smoothly on $M$ from the right. The action groupoid
$M\times G\rightrightarrows M$, denoted by $M\rtimes G$, is defined by the following data:
\begin{itemize}
\item $s(x,g)=x$;
\item $t(x,g)=xg$, so that
a pair $\Big((x,g),(x',g')\Big)$ is decomposable iff $x'=xg$;
\item $m\Big((x,g),(xg,g')\Big)=(x,gg')$;
\item  $i(x,g)=(xg,g^{-1})$;
\item $e(x)=(x,\textbf{1}_G)$.
\end{itemize}
\end{ex}

\subsection{Gerbes and $S^1$-central extensions of Lie groupoids}
\begin{ex}
Let $G$ be a Lie group and $BG$ its classifying stack. As we have seen, this is a stack, but it is in fact a rather special stack. This is because
\begin{enumerate}
\item Every manifold $X$  has at least one principal $G$ bundle over it, namely the trivial $G$ bundle;
\item Any two principal $G$ bundles are locally isomorphic.
\end{enumerate}
These two facts lead to the definition of a gerbe over a stack.
\end{ex}

\begin{defn}\label{gerbe2}
Let $\stack$ and $\gerbe$ be stacks over $\Sc$ and $\pi:\gerbe\To\stack$ a morphism of stacks.
Then $\pi:\gerbe\To\stack$ is called a \emph{gerbe} over (the stack) $\stack$, if
\begin{enumerate}
\item $\pi$ has local sections, i.e. there is an atlas $p: X\To\stack$ and a section $s:X\To\gerbe$ of $\pi\vert_X$,
where by a section we mean there exists a natural isomorphism
$\phi:\pi\circ s\Rightarrow p$ of functors.
\item Locally over $\stack$ all objects of $\gerbe$ are isomorphic, i.e. for any two objects $t_1,t_2\in\stack_T$ and
lifts $s_1,s_2\in\gerbe_T$ with $\pi(s_i)\isom t_i$, there is a covering $\{T_i\To T\}$ such that $s_1\vert_{T_i}\isom s_2\vert_{T_i}$.

\end{enumerate}
\end{defn}

A gerbe $\pi:\gerbe\To\stack$ is \emph{trivial}, if it admits a global section, i.e. if there exists a morphism of stacks $\sigma:\stack\To\gerbe$ satisfying
$\pi\circ\sigma\isom\id_\stack$.

\begin{defn}
A gerbe $\gerbe\To\stack$ is called an $S^1$-gerbe if there is an atlas $p:X\To\stack$ and a section $s:X\To\gerbe$ such that
there is an isomorphism
$$
\Phi:\Aut(s/p):=(X\times_\gerbe X)\times_{X\times_\stack X} X\isom S^1\times X
$$
as a family of groups over $X$ such that on $X\times_\stack X$ the diagram
$$
\xymatrix{
\Aut(s\circ\pr_1/p\circ\pr_1)\ar[rr]^\isom\ar[rd]^{\pr_1^\ast\Phi} & & \Aut(s\circ\pr_2/p\circ\pr_2)\ar[ld]_{\pr_2^\ast\Phi}\\
&X\times_\stack X\times S^1
}
$$
where the horizontal map is the isomorphism given by the universal property of the fibre product, commutes. This means that
the automorphism groups of objects of $\gerbe$ are central extensions of those of $\stack$ by $S^1$.
\end{defn}
\begin{defn}
Let $\Gamma=X_1 \rightrightarrows X_0$ be a Lie groupoid. An $S^1$-\emph{central extension} of $X_1\arr X_0$ consists of
\begin{enumerate}
\item a Lie groupoid $R_1\arr X_0$ and a morphism of Lie groupoids $(\pi,\id):[R_1\arr X_0]\To[X_1\arr X_0]$,
\item a left $S^1$ action on $R_1$, making $\pi:R_1\To X_1$ a left principal $S^1$ bundle. The action must satisfy
$(s\cdot x)(t\cdot y)=st\cdot(xy)$, for all $s,t\in S^1$ and $(x,y)\in R_1\times_{X_0}R_1$.
\end{enumerate}
\end{defn}

When $R_1\To X_1$ is topologically trivial, then $R_1\isom X_1\times S^1$ and the central extension is determined by a
\emph{groupoid $2$-cocycle} of $X_1\rightrightarrows X_0$ with values in $S^1$. This is a smooth map
$$
c: \Gamma^{(2)}=\Big\{(x,y)\in X_1\times X_1\mid s(x)=t(y)\Big\}\To S^1
$$
satisfying the cocycle condition
$$
c(x,y)c(xy,z)c(x,yz)^{-1}c(y,z)^{-1}=1
$$
for all $(x,y,z)\in\Gamma^{(3)}$. The groupoid structure on $R_1\rightrightarrows X_0$ is given by
$$
(x,\lambda_1)\cdot(y,\lambda_2)=(xy,\lambda_1\lambda_2 c(x,y)),
$$
for all $(x,y)\in\Gamma^{(2)}$ and $\lambda_1,\lambda_2\in S^1$.
\begin{prop}[Behrend, Xu, \cite{BXu}]\label{groupoidstogerbes}
Let $X_1\arr X_0$ be a Lie groupoid and $\stack$ its corresponding differential stack of $X_\bullet$-torsors. There is one-to-one correspondence
between $S^1$-central extensions of $X_1\arr X_0$ and $S^1$-gerbes $\gerbe$ over $\stack$ whose restriction to $X_0:\gerbe\vert_{X_0}$ admits a
trivialization.
\end{prop}

\subsection{Sheaf cohomology on differentiable stacks}

Let $\pi:\stack\To\Sc$ be a differentiable stack. 
Following \cite{La} and \cite{He} one can define sheaves of Abelian groups on $\stack$.

\begin{defn}
A sheaf $\F$ of Abelian groups on $\pi:\stack\To\Sc$ is
determined by the following data
\begin{enumerate}
\item For each morphism of stacks $X\To\stack$ where $X\in\Ob(\Sc)$ is a manifold, a sheaf $\F_{X\To\stack}$ of Abelian groups on $X$ in the usual
sense, i.e. an Abelian group $\F_{X\To\stack}(U)$ associated to each open $U\subseteq X$, etc.
\item For any $2$-commuting triangle
\begin{equation}
\xymatrix{
X\ar[rr]^f\ar[dr]_h &\ar@{}[d] |{\stackrel{\varphi}{\Longrightarrow}} & Y\ar[dl]^g\\
& \stack &
}
\end{equation}
with an isomorphism $\varphi:g\circ f\To h$ of functors, there exists a morphism of sheaves $\Phi_{f,\varphi}:f^\ast\F_{Y\To\stack}\To\F_{X\To\stack}$
(often denoted simply by $\Phi_f$) compatible for $X\To Y\To Z$. We require that $\Phi_f$ is an isomorphism, whenever $f$ is an open covering.
\end{enumerate}
The sheaf $\F$ is called \emph{Cartesian} if all $\Phi_f$ are isomorphisms.
\end{defn}

We
denote the category of Abelian sheaves on $\stack$ by $\Ab(\stack)$.

\begin{prop}
The category $\Ab(\stack)$ is an Abelian category with enough injective objects, i.e. for every object $\F\in\Ob(\Ab(\stack))$ there exists
an injection $0\To\F\To\I$ with $\I$ injective.
\end{prop}

\begin{defn}
Let $U$ be a manifold. A sheaf in the usual sense (i.e. defined only on open subsets of $U$) is called a \emph{small} sheaf on $U$.
\end{defn}
\begin{defn}
Let $\stack$ be a stack over $\Sc$ and $\F$ a sheaf over $\stack$. Let $x\in\Ob(\stack_U)$, where
$U\in\Ob(\Sc)$ is a manifold. The small sheaf on $U$, which maps the open subset $V\subseteq U$
to $\F(x\mid V)$ is called the small sheaf \emph{induced} by $\F$ via $x:U\To\stack$ on $U$. We denote it by $\F_{x,U}$ or simply
$\F_U$, if there is no risk of confusion.
\end{defn}

Given a morphism in $\theta:y\To x$ in $\stack$ lying over
a $\Cinfty$ map $f:V\To U$ in $\Sc$, there is an induced morphism of small sheaves over $V$
$$
\theta^\ast:f^{-1}\F_{x,U}\To\F_{y,V}.
$$
The cohomology of a sheaf $\F\in\Sh(\stack)$
is defined in the same way as it is defined for manifolds:
One first defines the \emph{global section} functor
$$
\Gamma(\stack,\cdot): \Ab(\stack)\To\Ab,
$$
where now
$$
\Gamma(\stack,\F):= \plim \Gamma(X,\F_{X\To\stack})
$$
and the limit is taken over all atlases $X\To\stack$, the transition functions for a $2$-commutative diagram
$\xymatrix{
X'\ar[rr]^f\ar[dr]_h &\ar@{}[d] |{\stackrel{\varphi}{\Longrightarrow}} & X\ar[dl]^g\\
& \stack &
}$
are
given by the restriction maps $\Phi_{f,\varphi}$.

Next one chooses an injective resolution $0\To\F\stackrel{\varepsilon}{\To}\I^\bullet$ and
sets
$$
H^i(\stack,\F)=h^i(\Gamma(\stack,\I^\bullet)).
$$
\begin{rem}
For a Cartesian sheaf $\F$ over $\stack$ the global section functor can be defined by choosing an atlas $X\To\stack$ and then setting
$$
\Gamma(\stack,\F):=\ker\Big(\Gamma(X,\F)\rightrightarrows\Gamma(X\times_\stack X)\Big).
$$
This is known to be independent of the chosen atlas $X\To\stack$ and moreover it coincides with the previous definition, \cite{He}.
\end{rem}

\begin{thm}[Giraud]\label{gerbeclass}
Isomorphism classes of $S^1$-gerbes over $\stack$ are in one-to-one correspondence with $H^2(\stack,\underline{S}^1)$.
\end{thm}

\subsection{\v Cech and simplicial cohomology of stacks}

\begin{defn}
Let $\Delta$ be the category whose objects are finite ordered sets $[n]=\{0<1<\cdots<n\}$, and whose morphisms
are nondecreasing monotone functions.
\end{defn}

\begin{defn}
Let $\mathcal{A}$ be a category. A \emph{simplicial object} $A$ in $\mathcal{A}$ is a contravariant functor
$A:\Delta^\op\To\A$
\end{defn}

\begin{defn}
A morphism of simplicial objects is a natural transformation between the corresponding functors, and the category
$\mathcal{SA}$ of all simplicial objects in $\A$ is just the functor category $\A^{\Delta^\op}$.
\end{defn}

\begin{prop}
To give a simplicial object $A$ in a category $\mathcal{A}$, it is necessary and sufficient to give a sequence of objects
$A_0,A_1,A_2,\ldots$ together with \emph{face operators} $\fa_i:A_p\To A_{p-1}$ and \emph{degeneracy operators}
$\de_i:A_p\To A_{p+1}$, where $i=0,1,\ldots, p$, satisfying the so called \emph{simplicial identities}:
\begin{eqnarray}
\fa_i\fa_j &=& \fa_{j-1}\fa_i,\qquad\textnormal{if } i<j\nonumber\\
\de_i\de_j&=&\de_{j+1}\de_i,\qquad\textnormal{if } i\leq j\nonumber\\
\fa_i\de_j &=& \left\{
\begin{array}{ll}
\de_{j-1}\fa_i, & \textnormal{if } i < j\\
\id, & \textnormal{if } i=j\textnormal{ or } i=j+1\\
\de_j\fa_{i-1}, & \textnormal{if } i>j+1.
\end{array}
\right.\nonumber
\end{eqnarray}
\end{prop}

\begin{proof}
Omitted. See \cite{We}, Prop. 8.1.3.
\end{proof}
If one dualizes the concept of simplicial objects, one obtains cosimplicial objects and the following proposition:

\begin{prop}
To give a cosimplicial object $A$ in a category $\mathcal{A}$, it is necessary and sufficient to give a sequence of objects
$A^0,A^1,\ldots$ together with \emph{coface operators} $\fa^i: A^{p-1}\To A^p$ and \emph{codegeneracy operators}
$\sigma^i: A^{p+1}\To A^p$,
where
$i=0,1,\ldots,p$, which satisfy the \emph{cosimplicial identities}
\begin{eqnarray}
\fa^j\fa^i &=& \fa^i\fa^{j-1},\qquad\textnormal{if } i<j\nonumber\\
\de^j\de^i &=& \de^i\de^{j+1},\qquad\textnormal{if } i\leq j\nonumber\\
\de^j\fa^i &=& \left\{
\begin{array}{ll}
\fa^i\de^{j-1}, & \textnormal{if } i < j\\
\id, & \textnormal{if } i=j\textnormal{ or } i=j+1\\
\fa^{i-1}\de^j, & \textnormal{if } i>j+1.
\end{array}
\right.\nonumber
\end{eqnarray}
\end{prop}

\begin{proof}
Omitted. See \cite{We}, Cor. 8.1.4.
\end{proof}

\begin{rem}\label{simfun}
It is clear by the above, that if we have a contravariant funtor $F:\mathcal{A}\To\mathcal{B}$, then $F$ maps simplicial objects in
$\mathcal{A}$ to cosimplicial objects in $\mathcal{B}$. In the same way, a covariant functor $F$ maps simplicial objects to simplicial objects, etc.
\end{rem}
\begin{defn}
Let $A$ be a simplicial object in an \emph{Abelian} category $\mathcal{A}$. The \emph{associated}, or \emph{unnormalized, chain complex}
$C(A)$ has its objects $C_p=A_p$, and its boundary morphism $d:C_p\To C_{p-1}$ is the alternating sum of the face operators
$\fa_i:C_p\To C_{p-1}$:
$$
d=\fa_0-\fa_1+\cdots+(-1)^p\fa_p.
$$
The simplicial identities for $\fa_i\fa_j$ imply that $d^2=0$, so that we indeed have a complex.
\end{defn}

We now come back to our original situation and define for all $p\geq 0$
$$
X_p=\underbrace{X\times_\stack\ldots\times_\stack X}_{p+1 \textnormal{ times}}.
$$
Since $X\To\stack$ is a representable submersion, all $X_p$ are manifolds. We want to make $X_\bullet=\{X_p\}$ into a
simplicial manifold, i.e. a simplicial object in the category of manifolds:

\begin{equation}\label{simpman}
\xymatrix{
\cdots
\ar@<2pt>[r]
\ar@<6pt>[r]
\ar@<-2pt>[r]
\ar@<-6pt>[r]
& X_2
\ar[r]
\ar@<4pt>[r]
\ar@<-4pt>[r]
& X_1
\ar@<2pt>[r]
\ar@<-2pt>[r] &
X_0.
}
\end{equation}
First, note that $X_p$ corresponds to the space of chains of composable $p$ arrows in the groupoid
$X_1\rightrightarrows X_0$.
Define the face and degeneracy maps so that
$$
\fa_i(g_1,\ldots,g_p) = \left\{
\begin{array}{ll}
(g_2,\ldots,g_p), & \textnormal{if } i =0\\
(g_1,\ldots,g_i g_{i+1},\ldots,g_n), & \textnormal{if } 0<i<p\\
(g_1,\ldots,g_{p-1}), & \textnormal{if } i=p,
\end{array}
\right.\nonumber
$$

$$
\sigma_i(g_1,\ldots, g_p)=(g_1,\ldots,g_i,1,g_{i+1},\ldots,g_p).
$$

\begin{ex}\label{Xpquotient}
We claim that for a quotient stack $[X/G]$ with the natural atlas $X\To[X/G]$
$$
X_p=\underbrace{X\times_\stack\ldots\times_\stack X}_{p+1 \textnormal{ times}}\isom X\times\prod_{i=1}^p G.
$$
This can be seen as follows. By definition $X_0=X$ and the product on the right hand side is empty, thus the claim is true when $p=0$.
Next note that by \cite{He} we have $X\times_\stack X\isom X\times G$. This implies that
\begin{eqnarray}
X\times_\stack X\times_\stack X&\isom& (X\times_\stack X)\times_X(X\times_\stack X)\isom(X\times G)\times_X(X\times G)\nonumber\\
&\isom& X\times G\times G.\nonumber
\end{eqnarray}
Here the last isomorphism follows since
$$
(X\times G)\times _X(X\times G)=\Big\{\Big((x_1,g_1),(x_2,g_2)\Big)\in (X\times G)\times(X\times G)\,\Big\vert\,\, x_1=x_2\Big\}.
$$
More generally, one may write for $p>2$
\begin{eqnarray}
X_{p+1}&=&\underbrace{X\times_\stack\ldots\times_\stack X}_{p+2 \textnormal{ times}}\isom
\underbrace{\Big(X\times_\stack\ldots\times_\stack X\Big)}_{p+1 \textnormal{ times}}
\times_X\Big(X\times_\stack X\Big)\nonumber\\
&\isom& X_p\times_X(X\times G)\isom X_p\times G\nonumber
\end{eqnarray}
and the claim follows from this by induction.
\end{ex}

Now, let $\F$ be a sheaf of Abelian groups on $\stack$. Every $X_p$ has $p+1$ canonical projections $X_p\To\stack$, which are all canonically
isomorphic to each other. We choose one of them and call it $\pi_p:X_p\To\stack$. Recall that $\pi_p$ as a map from a manifold to a stack can be
identified with an object of $\stack$ lying over $X_p$. We denote the Abelian group $\F(\pi_p)$ associated to the object $\pi_p$ by the
contravariant sheaf functor $\F$ by $\F(X_p)$. By Remark \ref{simfun} we have then a cosimplicial  Abelian group
\begin{equation}
\xymatrix{
\F(X_0)\ar@<2pt>[r]
\ar@<-2pt>[r] &
\F(X_1)\ar[r]
\ar@<4pt>[r]
\ar@<-4pt>[r] &
\F(X_2)\ar@<2pt>[r]
\ar@<6pt>[r]
\ar@<-2pt>[r]
\ar@<-6pt>[r]&
\cdots.
}
\end{equation}
Since the category of Abelian groups is an Abelian category, we may form the associated cochain complex
to $\F(X_\bullet)$:
\begin{equation}\label{homsim}
\xymatrix{
C(\F(X_\bullet)): & \F(X_0)\ar[r]^\fa & \F(X_1)\ar[r]^\fa & \F(X_2)\ar[r]^\fa & \cdots
}
\end{equation}

\begin{defn}
The homology groups of the complex (\ref{homsim}) are denoted by
$$
\check{H}^i(X_\bullet,\F)=h^i(\F(X_\bullet))
$$
and called the \emph{\v Cech} cohomology groups of $F$ with respect to the covering $X\To\stack$.
\end{defn}

As usual, there exists also a map $\check{H}^i(X_\bullet,\F)\To H^i(\stack,\F)$. Moreover, we have the following proposition

\begin{prop}
Let $\F$ be a Cartesian sheaf of Abelian groups on a differentiable stack $\stack$. Let $X\To\stack$ be an
atlas and $\F^\bullet$ the induced simplicial sheaf on the simplicial manifold $X_\bullet$. Then there is an
$E_1$-spectral sequence:
$$
E_1^{p,q}=H^q(X_p,\F_p)\Longrightarrow H^{p+q}(\stack,\F).
$$
Moreover,
$$
H^i(\stack,\F)\isom H^i(X_\bullet,\F^\bullet)
$$
for all $i\geq 0$, where the latter cohomology group is the simplicial cohomology of $\F^\bullet$.
\end{prop}

\begin{proof}
See \cite{De1}, \cite{He}.
\end{proof}

\begin{cor}\label{stacksim}
Let $\stack$ be a differentiable stack with an atlas $X\To\stack$. Then
$$
H^i(\stack,\underline{S}^1)\isom H^i(X_\bullet,\underline{S}^1)
$$
for all $i\geq 0$.
\end{cor}

\begin{ex}\label{stackclass}
Let again $\stack=[X/G]$ be the quotient stack and $\F=\underline{S}^1_\stack$. By Example \ref{Xpquotient}
$X_p\isom X\times\prod_{i=1}^p G$. Hence for each $p\geq 0$ the induced small sheaves of $\underline{S}^1$ on $X_p$ are
the sheaves $\underline{S}^1_{,X\times G^p}$. It follows now easily from Corollary \ref{stacksim} and \cite{Br1}, \cite{De1}, \cite{Go} that
the cohomology groups $H^i([X/G],\underline{S}^1)$ are isomorphic to the $G$-equivariant cohomology groups
of $X$. Especially, the group
$$
H^2([X/G],\underline{S}^1)\isom H^2(X\times G^\bullet,S^1_{X\times G^\bullet})
$$
classifies the isomorphism classes of $G$-equivariant gerbes on $X$ in the sense of Brylinski, \cite{Br1}.
\end{ex}

\subsection{Faddeev-Mickelsson anomaly in terms of differentiable gerbes and Lie groupoids}
This section contains our main results.

\subsubsection{Infinite-dimensional Lie groups of Mickelsson-Rajeev type}

\begin{defn}
Let $G$ be an I.L.H. (resp. Banach) Lie group (see Appendix A). An extension of $G$ by an I.L.H. (resp. Banach) Lie group $N$ is a short exact sequence with smooth
homomorphisms
$$
\xymatrix{
1\ar[r] &N\ar[r]^{i} &\hat{G}\ar[r]^{q} & G\ar[r] & 1
}
$$
and with a smooth local section $\sigma$ in the sense that there exists an open identity neighborhood $U\subseteq G$
on which $\sigma:U\To\hat{G}$ is smooth and $q\circ\sigma=\id_U$.
\end{defn}

\begin{rem}
One can use other classes of infinite dimensional manifolds and Lie groups in the definition as well,  see \cite{MK}.
\end{rem}

The infinite-dimensional Lie groups that we are interested in are those that appear in Yang-Mills theories as gauge transformation groups or
their extensions (\cite{PS}), \cite{Mi1} and \cite{AM}).

Let $\Hilb$ be a complex infinite dimensional separable Hilbert space with a given polarization $\Hilb=\Hilb_+\oplus\Hilb_-$, where
$\Hilb_{\pm}$ are closed subspaces of $\Hilb$. Let $\epsilon$ be the associated sign operator $\epsilon:\Hilb\To\Hilb,\,\epsilon^2=1$ and
$\epsilon\vert_{\Hilb_{\pm}}=\pm 1_{\Hilb_{\pm}}$. Let
$\GL(\Hilb)$ be the general linear group of $\Hilb$ consisting of all invertible  bounded linear operators of $\Hilb$.
\begin{defn}
We say that an infinite dimensional Lie group $\Ga$ is of \emph{Mickelsson-Rajeev type}, if it is of the form
$$
\Ga=\GL_{\mathcal{I}^p}:=\Big\{g=\Mat{a}{b}{c}{d}\in\GL(\Hilb)\mid [\epsilon,g]\in\mathcal{I}^{2p}\Big\},
$$
where $\mathcal{I}^p\subseteq\mathcal{K}(\Hilb)$ is a two-sided ideal in the algebra $\mathcal{B}(\Hilb)$, $p\in\N_+$, equipped with a Banach space topology 
$(\mathcal{I}^p,\norm{\cdot}_{\mathcal{I}^p})$ and
$\mathcal{I}^p\subseteq\mathcal{I}^q$ is dense in $\mathcal{I}^q$ whenever $p<q$. We define $\GL_{\mathcal{I}^p}$ to be a Banach-Lie group with topology given by the norm
$$
\norm{a}+\norm{b}_{\mathcal{I}^{2p}}+\norm{c}_{\mathcal{I}^{2p}}+\norm{d}.
$$
\end{defn}

We may extend the definition to the value $p=\infty$ by defining $\mathcal{I}^\infty:=\mathcal{K}(\Hilb)\subseteq\mathcal{B}(\Hilb)$. Then we have a
sequence of Banach-Lie groups
$$
\GL_{\mathcal{I}^1}\subseteq\GL_{\mathcal{I}^2}\subseteq\cdots\subseteq\GL_{\mathcal{I}^\infty}.
$$

\begin{ex}
One could choose for the $\mathcal{I}^p$'s the Schatten ideals $\Lp^p$ or the weak-$\Lp^p$ spaces $\Lp^{p+}$. 
\end{ex}

Let $\A$ be a contractible Banach manifold. We assume that there exists a set of maps $\Map(\A,S^1)$ such that this
set has a structure of a Banach-Lie group (compare with \cite{MiRa}, Remark on page 388).

We assume that our Lie group extension is of the form
$$
\hat{\Ga}=\widehat{\GL}_{\mathcal{I}^p}=(\Ev_p\times\Map(\A,S^1))/N,
$$
where
$$
\Ev_p=:\{(g,q)\mid g\in\GL_{\mathcal{I}^p},\, q\in\GL(\Hilb_+),\,aq^{-1}-1\in\mathcal{I}^p\}\subseteq\GL_{\mathcal{I}^p}\times\GL(\Hilb_+),
$$
$g=\Mat{a}{b}{c}{d}$ and the group multiplication is given by
$$
(g_1,q_1)(g_2,q_2)=(g_1g_2,q_1q_2).
$$
The topology of $\Ev_p$ is not the product space topology, but given by the norm
$$
\norm{(g,q)}=\norm{a}+\norm{d}+\norm{b}_{2p}+\norm{c}_{2p}+\norm{a-q}_p.
$$
Then $\Ev_p$ is a Banach-Lie group.
Above, $N$ is assumed to be a (closed) normal Banach-Lie subgroup
of $\Ev_p\times\Map(\A,S^1)$ consisting of elements of the
form $(1,q,\mu_q)$, where $\mu_q\in\Map(\A,S^1)$ depends smoothly on $q\in\GL(\Hilb_+)$. This makes $\hat{\Ga}$ into a Banach-Lie group.

The group $\widehat{\GL}_{\mathcal{I}^p}$ is assumed to be a (nontrivial) Banach principal $\Map(\A,S^1)$-bundle over $\GL_{\mathcal{I}^p}$ with the obvious projection map.
Near the unit element $1\in\GL_{\mathcal{I}^p}$ the formula
$$
\psi(g)=(g,a,1)\mod N,
$$
where $g=\Mat{a}{b}{c}{d}\in\GL_{\mathcal{I}^p}$, defines a local section $\psi:U\To\widehat{\GL}_{\mathcal{I}^p}$
of the principal $\Map(\A,S^1)$-bundle $p:\widehat{\GL}_{\mathcal{I}^p}\To\GL_{\mathcal{I}^p}$. 

\begin{defn}
An extension of infinite dimensional Lie groups $p:\hat{\Ga}\To\Ga$
is said to be of \emph{Mickelsson-Rajeev type} if it is of the above form.
\end{defn}

A Lie group extension 
of Mickelsson-Rajeev type defines a local $\Map(\A,S^1)$-valued (smooth) Lie group
2-cocycle $\omega$ by
$$
\psi(g_1)\psi(g_2)=\psi(g_1 g_2)(1,1,\omega(g_1,g_2)),
$$
where $\omega(g_1,g_2)\in\Map(\A,S^1)$. This can then be extended to a global $\Map(\A,S^1)$-valued (smooth) $2$-cocycle by translation giving an
element in the Lie group cohomology $[\omega]\in H^2(\GL_{\mathcal{I}^p},\Map(\A,S^1))$.

It follows from the definition that
$$
\lie(\widehat{\GL}_{\mathcal{I}^p})=\lie(\GL_{\mathcal{I}^p})\oplus\Map(\A,S^1),
$$
where the commutator in $\lie(\widehat{\GL}_{\mathcal{I}^p})$ is given by
$$
[(X,\mu),(Y,\nu)]=([X,Y],X\cdot\nu-Y\cdot\mu+\eta(X,Y;\cdot)),
$$
where $\eta$ is a $\Map(\A,S^1)$-valued Lie algebra cocycle on $\lie(\GL_{\mathcal{I}^p})$ and
the Lie derivative of a function $\nu$ on $\A$ to the direction of the vector field $X$ defined by the $\Ga$ action on $\A$ is denoted by $X\cdot\nu$.
Then at least in principle, one can calculate the Lie algebra cocycle $\eta$ as follows:
Let $\exp(tX)$ and $\exp(tY)$ be two one-parameter subgroups on $\GL_{\mathcal{I}^p}$. Then
$$
\frac{\partial^2}{\partial t\partial s}\psi(e^{tX})\psi(e^{sY})\psi(e^{-tX})\psi(e^{-sY})\Big\vert_{t=s=0}=([X,Y],0,\eta(X,Y)).
$$
\subsubsection{From principal $\Map(\A,S^1)$-bundles over $\Ga$ to line bundles over $\A\times\Ga$}

Let $\A$ be a contractible Banach manifold with a smooth right action of a Lie group $\Ga$ of Mickelsson-Rajeev type. 
We assume that a Lie group extension $p:\hat{\Ga}\To\Ga$ 
of Mickelsson-Rajeev type is given:
$$
\xymatrix{
\Map(\A,S^1)\ar@{-}[r] & \hat{\Ga}\ar[d]^p\\
& \Ga
}
$$
Here $p:\hat{\Ga}\To\Ga$  is a principal $\Map(\A,S^1)$-bundle.

Now, choose an open cover $\{U_\alpha\}_{\alpha\in I}$ of $\Ga$ and local sections $\psi_\alpha: U_\alpha\To\hat{\Ga}$. Over
the intersections $U_\alpha\cap U_\beta$, we have transition functions $\phi_{\alpha\beta}:U_\alpha\cap U_\beta\To\Map(\A,S^1)$
satisfying
$$
\psi_\alpha(g)=\psi_\beta(g)\phi_{\beta\alpha}(g),
$$
for all $g\in U_\alpha\cap U_\beta$. We can use the transition functions $\phi_{\alpha\beta}$ to construct a \emph{line} bundle over the product $\A\times\Ga$ as follows.
Define functions $\tilde{\phi}_{\beta\alpha}:(U_\alpha\cap U_\beta)\times\Ga\to S^1$ so that
$$
\tilde{\phi}_{\beta\alpha}(A,g):=\Big(\phi_{\beta\alpha}(g)\Big)(A)\in S^1,
$$
for all $A\in\A$ and $g\in\Ga$. The functions $\tilde{\phi}_{\beta\alpha}$ satisfy the following cocycle property
\begin{eqnarray}
\tilde{\phi}_{\gamma\beta}(A,g)\cdot\tilde{\phi}_{\beta\alpha}(A,g) &=&
\phi_{\gamma\beta}(g)(A)\cdot\phi_{\beta\alpha}(g)(A)\\
&=& \Big(\phi_{\gamma\beta}(g)\cdot\phi_{\beta\alpha}(g)\Big)(A)\nonumber\\
&=& \phi_{\gamma\alpha}(g)(A)\nonumber\\
&=& \tilde{\phi}_{\gamma\alpha}(A,g),\nonumber
\end{eqnarray}
and hence being transition functions determine an $S^1$-bundle over $\A\times\Ga$:
\begin{eqnarray}\label{line}
\xymatrix{
S^1\ar@{-}[r] & P\ar[d]^{\pi}\\
& \A\times\Ga
}
\end{eqnarray}
\begin{rem}
Note that the original $\Map(\A,S^1)$-bundle $p:\hat{\Ga}\To\Ga$ can be reconstructed from the transition functions of the
$S^1$-bundle $P\To\A\times\Ga$.
\end{rem}
\subsubsection{Constructing Lie groupoid operations on the line bundle over $\A\times\Ga$ -- The ``Cut and reglue'' procedure}
Suppose that the Mickelsson-Rajeev type Lie group extension
$p:\hat{\Ga}\To\Ga$
is given by the data of
a chosen open trivializing covering $\{U_\alpha\}$ of $\Ga$ with transition functions
$\phi_{\alpha\beta}:U_\alpha\cap U_\beta\To\Map(\A,S^1)$
and local 2-cocycles $\omega_{\alpha\beta,\gamma}:U_\alpha\times U_\beta\To\Map(\A,\Real)$ defining the multiplication on $\hat{\Ga}$
(this can always be done starting from the global extension and then looking at the trivializations). More precisely,
suppose that $f\in U_\alpha, g\in U_\beta$, $fg\in U_\gamma$ and $\lambda,\mu\in\Map(\A,S^1)$. Then the multiplication on the group $\hat{\Ga}$
is defined (locally) by the smooth maps
$$
m^{\hat{\Ga}}_{\alpha\beta,\gamma}:
\Big(U_\alpha\times\Map(\A,S^1)\Big)\times\Big(U_\beta\times\Map(\A,S^1)\Big)\To U_\gamma\times\Map(A,S^1),
$$
$$
m^{\hat{\Ga}}_{\alpha\beta,\gamma}\Big((f,\lambda),(g,\mu)\Big)=\Big(fg,\lambda(f\cdot\mu)\, e^{2\pi i \omega_{\alpha\beta,\gamma}(\cdot,f,g)}\Big),
$$
where $f\cdot\mu$ is the function $(f\cdot\mu)(A)=\mu(A^f)$ and for fixed $f$ and $g$
\begin{equation}\label{cocyc}
\omega_{\alpha\beta,\gamma}(\cdot;f,g):\A\To\Real,\quad\omega_{\alpha\beta,\gamma}(A;f,g):=\omega_{\alpha\beta,\gamma}(f,g)(A).
\end{equation}
Denoting $s_{\alpha\beta,\gamma}=e^{2\pi i\omega(\cdot,f,g)}$, the following  compatibility condition is satisfied:
\begin{equation}\label{gluemult}
s_{\alpha\beta,\gamma}(A;f,g)=\phi_{\alpha\alpha'}(A;f)\phi_{\beta\beta'}(A^f;g)\phi_{\gamma\gamma'}(A;fg)^{-1}s_{\alpha'\beta',\gamma'}(A;f,g),
\end{equation}
whenever  $f\in U_\alpha\cap U_{\alpha'},g\in U_\beta\cap U_{\beta'}$ and $fg\in U_\gamma\cap U_{\gamma'}$.
This is just the condition that we can glue together the local multiplication maps $m^{\hat{\Ga}}_{\alpha\beta,\gamma}$ to a
well-defined global smooth multiplication map
$m^{\hat{\Ga}}:\hat{\Ga}\times\hat{\Ga}\To\hat{\Ga}$.

Ignoring the various lower indices, the group 2-cocycle condition reads:
\begin{equation}\label{cocycle}
\omega(g_1g_2,g_3)+\omega(g_1,g_2)=\omega(g_1,g_2 g_3)+g_1\cdot\omega(g_2,g_3),
\end{equation}
where $g_1\cdot\omega(\cdot;g_2,g_3):\A\To\Real$ is the function
$$
g_1\cdot\omega(A;g_2,g_3)=\omega(A^{g_1};g_2,g_3).
$$
Notice, that this condition is equivalent to the associativity of the product on $\hat{\Ga}$.

Recall, that groupoid multiplication in
$\Gamma=(\A\rtimes\Ga\rightrightarrows\A;s,t,m,i,e)$ is defined by
\begin{eqnarray}
m:\Gamma^{(2)}&=& (\A\times\Ga)\times_{s,t} (\A\times\Ga)\To \A\times\Ga\nonumber\\
&=&\Big\{\Big((A_1,g_1),(A_2,g_2)\Big)\in(\A\times\Ga)\times(\A\times\Ga)\mid A_2=A_1^{g_1}\Big\}\To\A\times\Ga,\nonumber
\end{eqnarray}
$$
m\Big((A_1,g_1),(A_1^{g_1},g_2)\Big)=(A_1,g_1g_2),
$$
where
$$
s:\A\times\Ga\To\A,\quad s(A,g)=A
$$
is the source map and
$$
t:\A\times\Ga\To\A,\quad t(A,g)=A^g.
$$
is the target map.

Now $\{\A\times U_\alpha\}_{\alpha\in I}$ is an open covering of $\A\times\Ga$. We use the local group 2-cocycles
$\omega_{\alpha\beta,\gamma}:U_\alpha\times U_\beta\To\Map(\A,\Real)$ to define maps
$$
c_{\alpha\beta,\gamma}:
\Big\{\Big((A_1,g_1),(A_2,g_2)\Big)\in(\A\times U_\alpha)\times(\A\times U_\beta)\,\Big\vert\, A_2=A_1^{g_1},\, g_1 g_2\in U_\gamma\Big\}\To S^1,
$$
$$
c_{\alpha\beta,\gamma}(A_1,g_1,A_1^{g_1},g_2)=e^{2\pi i \omega_{\alpha\beta,\gamma}(A_1,g_1,g_2)}.
$$
We assume that the 2-cocycles $\omega_{\alpha\beta,\gamma}$ depend smoothly on the variable $A\in\A$ so that the maps $c_{\alpha\beta,\gamma}$
are smooth as well, when we give the sets where the different $c_{\alpha\beta,\gamma}$ are defined the manifold structure described below.
It follows from (\ref{cocycle}) that these satisfy the following cocycle condition
\begin{eqnarray}\label{ass}
& &c(A_1,g_1,A_1^{g_1},g_2)c(A_1,g_1 g_2, A_1^{g_1g_2},g_3)=c(A_1^{g_1},g_2,A_2^{g_2},g_3)\cdot\\
& &c(A_1,g_1,A_1^{g_1},g_2 g_3).\nonumber
\end{eqnarray}

Next, we define the following local multiplication maps by
\begin{eqnarray}
m_{\alpha\beta,\gamma}:
& &\Big\{\Big((A_1,g_1,\lambda),(A_2,g_2,\mu)\Big)\in(\A\times U_\alpha\times S^1)\times(\A\times U_\beta\times S^1)\,\Big\vert\, A_2=A_1^{g_1},\nonumber\\
& &g_1\in U_\alpha,\, g_2\in U_\beta,\, g_1 g_2\in U_\gamma,\, \lambda,\mu\in S^1\Big\}\To \A\times U_\gamma\times S^1,\nonumber
\end{eqnarray}
$$
m_{\alpha\beta,\gamma}\Big((A_1,g_1,\lambda),(A_1^{g_1},g_2,\mu)\Big)=\Big(A_1,g_1g_2,\lambda\mu\cdot c_{\alpha\beta,\gamma}(A_1,g_1,A_1^{g_1},g_2)\Big).
$$
Notice, that the set where $m_{\alpha\beta,\gamma}$ is defined is an open subset of the manifold
$$
(\A\times U_\alpha\times S^1)\times_{s\circ\pr_{1,2};\A;t\circ\pr_{1,2}}(\A\times U_\beta\times S^1)
$$
as the inverse image of the open set $U_\gamma\subseteq\Ga$ under the smooth map
$$
m^{\alpha\beta}:(\A\times U_\alpha\times S^1)\times_{s\circ\pr_{1,2};\A;t\circ\pr_{1,2}}(\A\times U_\beta\times S^1)\To\Ga,
$$
$$
m^{\alpha\beta}\Big((A_1,g_1,\lambda),(A_1^{g_1},g_2,\mu)\Big)=g_1g_2.
$$
Moreover, $(\A\times U_\alpha\times S^1)\times_{s\circ\pr_{1,2};\A;t\circ\pr_{1,2}}(\A\times U_\beta\times S^1)$ is indeed a manifold, since
both maps $s\vert_{\A\times U_\alpha}\circ\pr_{1,2}$ and $t\vert_{\A\times U_\beta}\circ\pr_{1,2}$ are surjective submersion as composites of surjective submersions.
Similarly, each $c_{\alpha\beta,\gamma}$ is defined on an open subset of the manifold
$$
(\A\times U_\alpha)\times_{s\vert{\A\times U_\alpha};\A;t\vert_{\A\times U_\beta}}(\A\times U_\beta)
$$

Since the restrictions $P\vert_{\A\times U_\alpha}=:\pi^{-1}(\A\times U_\alpha)\To\A\times U_\alpha$ of the $S^1$-bundle $P$ in (\ref{line}) are trivial, i.e.
there exists an $S^1$-bundle isomorphism
$$
P\vert_{\A\times U_\alpha}\isom \A\times U_\alpha\times S^1,
$$
one can patch together the various maps $m_{\alpha\beta,\gamma}$ to obtain a partial multiplication map $m_P$ on
the total space $P$ of the $S^1$-bundle $\pi:P\To\A\times\Ga$. Here by ``partial multiplication'' we mean that not every pair of elements in $P$ can
be multiplied together. The cocycle condtion (\ref{ass}) guarantees that the multiplication map $m_P$ is associative.
We want to make these arguments rigorous and show, that this makes $P\rightrightarrows\A$ a groupoid.
\begin{prop}
$(P\rightrightarrows\A,m_P,s_P,t_P)$ is a Banach-Lie groupoid, where the source and target map $s_P$ and $t_P$ are defined so that
$$
s_P=s\circ\pi\qquad t_P=t\circ\pi.
$$
\end{prop}
\begin{proof}
First, note that $s_P$ and $t_P$ are surjective submersions as compositions of two surjective submsersions.

Next, choose bundle isomorphisms giving local trivializations
$$
\varphi_\alpha:\A\times U_\alpha\times S^1\stackrel{\sim}{\To} P\vert_{\A\times U_\alpha},
$$
for each $\alpha\in I$. Hence for each $\alpha\in I$ we have a commutative diagram
$$
\xymatrix{
\A\times U_\alpha\times S^1\ar[r]_{\isom}^{\varphi_\alpha}\ar[d]_{\pr_{1,2}} & P\vert_{\A\times U_\alpha}\ar[d]^{\pi\vert_{\A\times U_\alpha}}\\
\A\times U_\alpha\ar[r]^{\id} & \A\times U_\alpha
}
$$
where $\varphi_\alpha$ is an $S^1$-equivariant map of manifolds and $\pr_{1,2}(A,g,\lambda)=(A,g)$. From this we see that
$$
s_P\vert_{\A\times U_\alpha}\circ\varphi_\alpha=s\vert_{\A\times U_\alpha}\circ\pr_{1,2}=\pr_1,
$$
where
$$\pr_1:\A\times U_\alpha\times S^1\To\A,\quad\pr_1(A,g,\lambda)=A.
$$
and
$$
s_P\vert_{\A\times U_\alpha}:P\vert_{\A\times U_\alpha}\To\A,\quad s\vert_{\A\times U_\alpha}:\A\times U_\alpha\To\A,
$$
$$
s_P\vert_{\A\times U_\alpha}=s\vert_{\A\times U_\alpha}\circ \pi\vert_{\A\times U_\alpha}.
$$

Hence
$$
s_P\vert_{\A\times U_\alpha}=\pr_1\circ\varphi_\alpha^{-1}.
$$
Similarly
$$
t_P\vert_{\A\times U_\alpha}\circ\varphi_\alpha=t\vert_{\A\times U_\alpha}\circ\pr_{1,2}
$$
or
$$
t_P\vert_{\A\times U_\alpha}=t\vert_{\A\times U_\alpha}\circ\pr_{1,2}\circ\varphi_\alpha^{-1}.
$$

We want to construct a global multiplication map
$$
m_P:P\times_{s_P,\A,t_P} P\To P
$$
from the local multiplication maps $m_{\alpha\beta,\gamma}$ introduced above. We denote by
$s_{P,\alpha}=s_P\vert_{\A\times U_\alpha}$ for every $\alpha\in I$ and similarly $t_{P,\alpha}=t_P\vert_{\A\times U_\alpha}$. Then
$$
P\vert_{\A\times U_\alpha}\times_{s_{P,\alpha;\A;t_{P,\beta}}} P\vert_{\A\times U_\beta}\subseteq P\times_{s_P,\A,t_P} P.
$$
Define
$$
\Big(P\vert_{\A\times U_\alpha}\times_{s_{P,\alpha;\A;t_{P,\beta}}} P\vert_{\A\times U_\beta}\Big)_\gamma
$$
as the open subsetset of $P\vert_{\A\times U_\alpha}\times_{s_{P,\alpha;\A;t_{P,\beta}}} P\vert_{\A\times U_\beta}$ so that
$$
\Big(P\vert_{\A\times U_\alpha}\times_{s_{P,\alpha;\A;t_{P,\beta}}} P\vert_{\A\times U_\beta}\Big)_\gamma:=
\Big(m^{\alpha\beta}\circ(\pr_2\times\pr_2)\circ(\varphi_\alpha^{-1}\times\varphi_\beta^{-1})\Big)^{-1}(U_\gamma).
$$
We may now define $m_{P;\alpha\beta,\gamma}:\Big(P\vert_{\A\times U_\alpha}\times_{s_{P,\alpha;\A;t_{P,\beta}}} P\vert_{\A\times U_\beta}\Big)_\gamma
\To P,$
$$
m_{P;\alpha\beta,\gamma}=\varphi_\gamma\circ m_{\alpha\beta,\gamma}\circ(\varphi_\alpha^{-1}\times\varphi_\beta^{-1}).
$$
This gives us a well-defined global multiplication map $m_P:P\times_{s,t} P\To P$, because of equation (\ref{gluemult}), that guarantees us that
the local multiplication maps at the group extension level glue together.

The other maps in the definition of a Lie groupoid are defined on local
trivializations $P\vert_{\A\times U_\alpha}\isom \A\times
U_\alpha\times S^1$ so that
$$
e_P(A)=(A,1_G,1),
$$
$$
i_P(A,g,\lambda)=(A^g,g^{-1},\lambda^{-1})
$$
\end{proof}
\begin{prop}
$P\rightrightarrows\A$ is an $S^1$-(Banach-Lie groupoid) central extension of the action gropoid $\A\rtimes\Ga$.
\end{prop}
\begin{proof}
We first claim, that the following diagrams commute:
\begin{equation}\label{first}
\xymatrix{
P\ar@<2pt>[d]^{t_P}\ar@<-2pt>[d]_{s_P}\ar[r]^\pi & \A\times\Ga\ar@<2pt>[d]^{t}\ar@<-2pt>[d]_{s}\\
\A\ar[r]^{\id} & \A
}
\end{equation}
\begin{equation}\label{second}
\xymatrix{
P\ar[r]^\pi & \A\times\Ga\\
\A\ar[r]^\id\ar[u]^{e_P} & \A\ar[u]_e
}
\end{equation}
\begin{equation}\label{third}
\xymatrix{
P\times_{s_P,t_P} P\ar[r]^{\pi\times\pi}\ar[d]_{m_P}& \,\,(\A\times\Ga)\times_{s,t}(\A\times\Ga)\ar[d]^m \\
P\ar[r]^\pi & \A\times\Ga
}
\end{equation}
\begin{equation}\label{fourth}
\xymatrix{
P\ar[d]_{i_P}\ar[r]^\pi &\A\times\Ga\ar[d]^i\\
P\ar[r]^\pi &\A\times\Ga
}
\end{equation}
\begin{enumerate}
\item Now, diagram (\ref{first}) commutes by definition.
\item On local trivializations of the $S^1$-bundle
$\pi:P\To\A\times\Ga$, the elements of the total space $P$ are of the form
$(A,g,\lambda)$, where $A\in\A$, $g\in\Ga$ and $\lambda\in S^1$.  Hence
$$
(\pi\circ e_P)(A)=\pi(A,1_G,1)=(A,1_G)=e(A),
$$
so that (\ref{second}) commutes.
\item Again, locally
$$
m_P\Big((A_1,g_1,\lambda),(A_1^{g_1},g_2,\mu)\Big)=\Big(A_1,g_1g_2, c(A_1,g_1,A_1^{g_1},g_2)\Big)
\stackrel{\pi}{\mapsto} (A_1,g_1g_2),
$$
and on the other hand
$$
(\pi\times\pi)\Big((A_1,g_1,\lambda),(A_1^{g_1},g_2,\mu)\Big)=\Big((A_1,g_1),(A_1^{g_1},g_2)\Big)
\stackrel{m}{\mapsto} (A_1,g_1g_2),
$$
which shows that (\ref{third}) commutes.
\item On local trivializations
$$
(i\circ\pi)(A,g,\lambda)=i(A,g)=(A^g,g^{-1})=\pi(A^g,g^{-1},\lambda^{-1})=(\pi\circ i_P)(A,g,\lambda).
$$
\end{enumerate}
This data gives us a morphism of Lie groupoids $(\pi,\id):[P\rightrightarrows\A]\To[\A\rtimes\Ga\rightrightarrows\A]$. Moreover,
$\pi:P\To\A\rtimes\Ga$ is a principal $S^1$-bundle by construction. The only thing left is to check that
$(s\cdot x)(t\cdot y)=(st)\cdot(xy)$ for all $s,t\in S^1$ and $(x,y)\in P\times_{s_P,\A,t_P} P$.
To see this, we look at the local picture, again. Thus, let $x=(A_1,g_1,\lambda)$ and $y=(A_1^{g_1},g_2,\mu)$. Now
\begin{eqnarray}
(s\cdot x)(t\cdot y) &=& (A_1,g_1,s\lambda)\cdot(A_1^{g_1},g_2,\mu)=\Big(A_1,g_1g_2,st\lambda\mu\cdot c(A_1,g_1,A_1^{g_1},g_2)\Big)\nonumber\\
&=&(st)\cdot(xy).\nonumber
\end{eqnarray}
\end{proof}
By Example 2.26. in \cite{Xu} the cocycle condition (\ref{ass}) of the family $\{c_{\alpha\beta,\gamma}\}$ guarantees that it
gives a $2$-cocycle in the simplicial cohomology $H^2(\A\times\Ga^\bullet,\underline{S}^1)$ (i.e. an element of the \v Cech cohomology with respect to out groupoid cover).
On the other hand this class is the class corresponding to the Morita equivalence class of the constructed $S^1$-groupoid extension
of $\A\rtimes\Ga$ under the isomorphism
$$
\Ext^{sm}(\A\rtimes\Ga,S^1)\isom H^2(\A\times\Ga^\bullet,\underline{S}^1).
$$
(see Proposition 2.17, \cite{Xu}). Next, recall from Example \ref{stackclass} that the Lie groupoid $\A\rtimes\Ga$ corresponds to the quotient stack $[\A/\Ga]$ and
$$
H^2(\A\times\Ga^\bullet,\underline{S}^1)\isom H^2([\A/\Ga],\underline{S}^1).
$$
Propostion \ref{groupoidstogerbes} produces then a gerbe $\gerbe$ over the stack $[\A/\Ga]$ whose gerbe class is the cohomology class of the
$2$-cocycle $\{c_{\alpha\beta,\gamma}\}$.

\begin{rem}
Note that the original multiplication in $\hat{\Ga}$ can be reconstructed from the associated $S^1$-groupoid extension $P\rightrightarrows\A$ using
(\ref{cocyc}). Since we noticed earlier that the original $\Map(\A,S^1)$-bundle $p:\hat{\Ga}\To\Ga$ can be reconstructed from the associated
$S^1$-bundle $\pi:P\To A\times\Ga$ we conclude that the whole group extension $\hat{\Ga}$ with its original principal bundle
structure can be reconstructed from the associated $S^1$-gropoid extension $P\rightrightarrows\A$.
\end{rem}

\appendix

\section{I.L.H. manifolds and Lie groups}

Our references are \cite{Br2} and \cite{Pa}.

\begin{defn}
A topological vector space $E$ is called an I.L.H. vector space if $E=\varprojlim_n \Hilb_n$ is an inverse limit
of separable Hilbert spaces $\Hilb_n$.
\end{defn}

Hence, the topology of an I.L.H. vector space $E$ is the inverse limit topology. This is the coarsest topology which makes
all the projection maps $p_n:E\To\Hilb_n$ continuous. Often one wants to impose the following extra condition in the definition of an
I.L.H. vector space:
\begin{itemize}
\item For every open ball $B$ in $\Hilb_n$, we have
\begin{equation}\label{ILH}
p_n^{-1}(\overline{B})=\overline{p_n^{-1}(B)}.
\end{equation}
\end{itemize}

\begin{thm}
Let $X$ be a paracompact manifold, modelled on an I.L.H. vector space $E$ satisfying (\ref{ILH}). Then for any
open covering $\U=\{U_i\}_{i\in I}$ of $X$ there exists a smooth partition of unity subordinate to $\U$.
\end{thm}

\begin{defn}
An I.L.H. topological group $G$ is called an \emph{I.L.H. Lie group} if it is a smooth I.L.H. manifold with the
group operations given by smooth I.L.H. maps.
\end{defn}

\begin{defn}
Let $P, B$ be smooth I.L.H. manifolds modelled on I.L.H. vector spaces $E$ and $F$ respectively,
$\pi:P\To B$ a smooth I.L.H. map and $G$ an I.L.H. Lie group. Then
$(P,B,G,\pi)$ is an I.L.H. principal bundle if the transition maps are smooth I.H.L. maps.
\end{defn}

Let $(P,M,G,\pi)$ be a smooth principal $G$-bundle on a closed manifold $M$, where we assume all the manifolds to be finite
dimensional and that $G$ is compact. Let $E=\ad P:=P\times_G\lie(G)$, where $G$ acts on $\lie(G)$ by the adjoint action, and $F:=T^\ast M\tensor\ad P$.

\begin{ex}
The space $\A(P)$ of smooth connections on $P$ is an affine I.L.H space with tangent vector space $\Cinfty(F)$.
\end{ex}

\begin{ex}
Let $E_G=\Ad P:=P\times_G G$ where $G$ acts on itself by the adjoint action. Then the set $\Ga(P):=\Cinfty(E_G)$ is an I.L.H. Lie group
modelled on $\Cinfty(E)$. It corresponds to the group of \emph{gauge transformations} of the principle $G$-bundle $P$, i.e. the group of
automorphisms of $P$ that cover the identity.
\end{ex}

\begin{ex}[Infinite dimensional Grassmannian of Segal and Wilson]
Let $\Hilb$ be a separable Hilbert space with an orthogonal decomposition $\Hilb=\Hilb_+\oplus\Hilb_-$. Recall that
for any two Hilbert spaces $\Hilb_1$ and $\Hilb_2$ the space $H.S.(\Hilb_1,\Hilb_2)$ of Hilbert-Schmidt operators $T:\Hilb_1\To\Hilb_2$ is
a Hilbert space with norm $\norm{T}_2=\sqrt{\tr(T^\ast T)}$. Let $\Gr_{res}(\Hilb)$ denote the set of closed
subspaces $W\subseteq\Hilb$ such that
\begin{enumerate}
\item The orthogonal projection onto $\Hilb_+,\,\pr_W^+:W\To\Hilb_+$ is Fredholm;
\item The orthogonal projection onto $\Hilb_-,\,\pr_W^-:W\To\Hilb_-$ is Hilbert-Schmidt.
\end{enumerate}
Then $\Gr_{res}(\Hilb)$ is a Hilbert
manifold modelled on $H.S.(\Hilb_+,\Hilb_-)$.
\end{ex}

\bibliographystyle{alpha-loc}

\end{document}